\documentclass[12pt]{amsart}
\usepackage{amssymb}
\usepackage{tikz}
\usepackage[symbol]{footmisc}

\usetikzlibrary{arrows, decorations.pathmorphing,
backgrounds,positioning,fit,petri}
\input amssym.def

 2

\newcommand{\R}{{\bf R}}
\newcommand{\bP}{{\bf P}}
\newcommand{\rR}{{\rm R}}
\newcommand{\C }{{\bf C}}
\newcommand{\rC }{{\rm C}}
\newcommand{\Z}{{\bf Z}}
\newcommand{\F}{{\rm F}}
\newcommand{\Q}{{\bf Q}}

\newcommand{\rN}{{\rm N}}
\newcommand{\A}{{\rm A}}
\newcommand{\D}{{\rm D}}

\newcommand{\SL}{{\bf SL}}
\newcommand{\E}{{\rm E}}
\newcommand{\mM}{{\mathcal{M}}}
\newcommand{\B}{{\rm B}}
\newcommand{\M}{{\rm M}}
\newcommand{\rP}{{\rm P}}
\newcommand{\T}{{\rm T}}
\newcommand{\K}{{\rm K}}
\newcommand{\I}{{\rm I}}

\renewcommand{\k}{{ k}}
\newcommand{\gH}{{\goth{H}}}
\newcommand{\G}{{\rm G}}
\newcommand{\X}{{\rm X}}
\newcommand{\Y}{{\rm Y}}
\newcommand{\bZ}{{\rm Z}}

\renewcommand{\Im}{{\rm Im}}
\renewcommand{\Re}{{\rm Re}}
\renewcommand{\ge}{{~\geqslant~}}
\renewcommand{\le}{{~\leqslant~}}

\def\bfminus{\mathrel{\hbox{\bf -\kern-1pt-}}}
\def\ckomma{\, \raise2pt \hbox{,}}

\begin{document}

\title[Critical points of Eisenstein series]{Critical points of 
Eisenstein series}

\author{Sanoli Gun and Joseph Oesterl\'e}

\address[Sanoli Gun]
      {The Institute of Mathematical Sciences,
        HBNI, C.I.T. Campus, Taramani, 
        Chennai 600113, India.}
\email[Sanoli Gun ]{sanoli@imsc.res.in}

\address[Joseph Oesterl\'e]
      {Institut de Math\'ematique de Jussieu - Paris Rive Gauche,      
       Sorbonne Universit\'e, Campus Jussieu, 
       4 place Jussieu,  75005-Paris, France.}
\email[J. Oesterl\'e]{joseph.oesterle@imj-prg.fr}


  \phantom{m}

\thanks{Research of this article was partially supported
by the Indo-French Program in Mathematics (IFPM).  
Both authors would like to thank IFPM for financial support
and Institute of Mathematical Sciences and Sorbonne Universit\'e
for providing excellent research environments.
The first author would also like to acknowledge MTR/2018/000201,
SPARC project~445 and DAE number theory plan project
for partial financial support.}

\maketitle

\begin{abstract}
For any even integer $k \ge 4$, let $\E_k$ be the normalized Eisenstein
series of weight $k$ for $\SL_2(\Z)$. Also let $\D$ be the 
closure of the standard fundamental domain of the Poincar\'e upper 
half plane modulo  $\SL_2(\Z)$. F.~K.~C.~Rankin and 
H. P. F. Swinnerton-Dyer showed that all zeros of 
$\E_k$ in $\D$ are of modulus one. In this article, we 
study the critical points of $\E_k$, that is to say the zeros of
the derivative of $\E_k$. We show that they are
simple. We count those belonging to $\D$, prove that they are located on the 
two vertical edges of $\D$ and produce explicit intervals that
separate them. We then count those belonging to $\gamma\D$,
for any $\gamma \in \SL_2(\Z)$.
\end{abstract}

\smallskip

\section{Introduction and statement of the results}

\medskip
\subsection{The Eisenstein series \nopunct}$~~$\\

Let $k$ be an even integer $\ge 4$. For any point $z$ in the Poincar\'e
upper half-plane $\mathfrak{H}$, let
\begin{equation}\label{Eis}
\G_{k}(z) := \sum_{ (c, d) \in \Z^2 \atop (c, d) \ne (0,0)} (cz + d)^{-k}.
\end{equation}
We have
\begin{equation}\label{new1}
\G_k(z) = 2 \zeta(k) \E_k(z),
\end{equation}
where $\E_k$ is {\it the normalized Eisenstein series of weight $k$
for} $\SL_2(\Z)$. Its Fourier expansion is given by
\begin{equation}\label{NEis}
\E_k(z) =  1 - \frac{2 k}{\B_k} \sum_{n=1}^{\infty} \sigma_{k-1}(n) q^n,
\end{equation}
where $q :=  e^{2\pi i z}$, $\sigma_{k-1}(n)$ is the sum of the 
$(k-1)$-th powers of the positive divisors of $n$ and $\B_k$ is the 
Bernoulli number of index $k$.

\smallskip
We recall in subsection 1.2 the results of Rankin and Swinnerton-Dyer
concerning the location of zeros of $\E_k$, or what amounts to the same
of $\G_k$. We then state in subsections 1.3 and 1.4, the 
results we obtain concerning the location of the critical points
of $\E_k$, i.e. of the zeros of the derivative $\E_{k}'$ of $\E_k$.

\medskip
\subsection{Results of F.K.C. Rankin and 
H.P.F. Swinnerton-Dyer \nopunct}$~~$\\

Since the function $\E_k$ is a modular form of weight $k$ for
the group $\SL_2(\Z)$, the set of its zeros in $\mathfrak{H}$ 
is stable under the action of $\SL_2(\Z)$, and the zeros of
a same orbit have the same multiplicity. We therefore content
ourselves with describing those belonging to a fundamental domain 
of $\mathfrak{H}$ modulo $\SL_2(\Z)$.  

\smallskip
Let us denote by $\D$ the closure in $\mathfrak{H}$ of the
standard fundamental domain of $\mathfrak{H}$
modulo $\SL_2(\Z)$, i.e. the set of points $z$ in $\mathfrak{H}$
such that $|\Re(z)| \le \frac{1}{2}$ and $|z| \ge 1$.
Rankin and Swinnerton-Dyer proved in \cite{RSD} that
all the zeros of $\E_k$ in $\D$ are of modulus $1$, and more 
precisely that : \par

\smallskip
$a)$ The points $e^{\pi i/3}$ and $e^{2\pi i / 3}$ are simple
zeros of $\E_k$ when $k\equiv 4\bmod 6$, 
zeros of multiplicity $2$ when $k \equiv 2 \bmod 6$ and are not 
zeros of $\E_k$ when $ k\equiv 0 \bmod 6$. \par

\smallskip
$b)$ The point $i = e^{\pi i/2}$ is a simple zero 
of $\E_k$ when $k \equiv 2 \bmod 4$ and 
is not a zero of $\E_k$
when $k \equiv 0 \bmod 4$. \par

\smallskip
$c)$ The zeros of $\E_k$ lying on the open circular 
arc consisting of the $e^{i \theta}$, with $\frac{\pi}{3} < \theta <~\frac{\pi}{2}$,
are simple. Their number is $[\frac{k}{12}]$ when
$k \not\equiv 2 \bmod 12$ and $[\frac{k}{12}] - 1$
when $k \equiv 2 \bmod 12$. Here $[x]$ 
denotes the largest integer less
than or equal to~$x$. The zeros located
on the open circular arc  consisting of the 
$e^{i \theta}$, with $\frac{\pi}{2} < \theta <~\frac{2\pi}{3}$,
are symmetric to the previous ones with respect to the
imaginary axis, and are simple as well.

\bigskip
\noindent
{\bf Remark 1.}$~-$ These results are consistent with the
following fact : since $\E_k$ is a modular form of weight
$k$ which does not vanish at infinity, the 
weighted number of its zeros in $\mathfrak{H}$ modulo 
$\SL_2(\Z)$ (counted with multiplicities, with weight~$\frac{1}{2}$ 
for those in the orbit of $e^{\pi i/2}$, $\frac{1}{3}$ for those in the
orbit of $e^{i\pi/3}$, and $1$ otherwise) has to be equal
to $\frac{k}{12}$.

\medskip
\noindent
{\bf Remark 2.}$~-$
It follows from \eqref{NEis} that 
the function $t \mapsto \E_k( \frac{1}{2} + it)$ is real-valued 
and analytic for $t >0$, and tends to
$1$ when $t$ tends to $+ \infty$. Further, it does 
not vanish for $t > \frac{\sqrt{3}}{2}$ by the results stated above.
Hence we have $\E_k (\frac{1}{2} + it) > 0$ for all $t > \frac{\sqrt{3}}{2}$.

\medskip
\subsection{Zeros of the derivative of $\E_k$ in $\D$ \nopunct}$~~$\\

In this article, we are interested in the zeros of the derivative of $\E_k$.  
This derivative $\E_k'$ is not a modular form, but only a quasi-modular form.
The set of its zeros is no longer stable under the action of $\SL_2(\Z)$.
Thus in order to get precise information about these zeros, we need to 
count them not only in $\D$, but also in $\gamma\D$ for any
$\gamma \in \SL_2(\Z)$. We state in this subsection three results 
about the zeros located in $\D$, and in the next subsection, those about 
the zeros located in $\gamma \D$, for any $\gamma \in \SL_2(\Z)$.

\bigskip
\noindent{\bf Theorem 1.}$-$
{\em The zeros of $\E_k'$ in $\D$ are simple, and 
have real part either $\frac{1}{2}$ or $-\frac{1}{2}$.
Those with real part $\frac{1}{2}$ are the translates
by $1$ of those with real part~$-\frac{1}{2}$.}

\bigskip
The second assertion follows from the fact that 
$\E_k'(z+1) = \E_k'(z)$.  So let us restrict our attention
to the zeros of $\E_k'$ in $\D$ with real part $\frac{1}{2}$,
i.e. located on the closed half-line 
$\frac{1}{2} + i [\frac{\sqrt{3}}{2}, +\infty[$.  Their
number is given by the following theorem :

\bigskip
\noindent{\bf Theorem 2.}$-$ 
{\em 
$a)$  The point $\frac{1}{2} + i \frac{\sqrt{3}}{2}$ 
is a zero of $\E_k'$ if and only if $k \equiv 2 \bmod 6$. \par
$b)$
The number of zeros of $\E_k'$ in the open
half-line $\frac{1}{2} + i ]\frac{\sqrt{3}}{2}, +\infty[$ 
is $[ \frac{k-4}{6} ]$.
}

\bigskip
Assertion $b)$ can be rephrased by saying that the number 
of zeros of $\E_k'$ in $\frac{1}{2} + i ]\frac{\sqrt{3}}{2}, +\infty[$ is 
$[\frac{k}{6}]$ when $k \equiv 4 \bmod 6$ 
and $[\frac{k}{6}] -1$ otherwise.  The following theorem
gives more precise information about their locations 
by producing open intervals that separate their imaginary 
parts :

\bigskip
\noindent{\bf Theorem 3.}$-$
{\em 
Let $\M := [\frac{k}{6}]$. For $1 \le m \le \M$,
let $t_m := \frac{1}{2}\cot(\frac{m \pi}{k+1})$.
We have therefore
$$
t_1 > \cdots > t_{\M} > \frac{\sqrt{3}}{2}.
$$
For each integer $m$ such that $1 \le m \le \M-1$, 
$\E_k'$  has a unique zero in the set
$\frac{1}{2}~+~i~]t_{m+1}, t_m[$.
Moreover, when $k \equiv 4 \bmod 6$ and $k \ne 4$,
$\E_k'$ has a unique zero in $\frac{1}{2} + i ]\frac{\sqrt{3}}{2}, t_{\M}[$.}

\bigskip
Let us briefly outline the plan of the proofs of the above
three theorems, which will be given in full detail in sections 2 
and 3. We shall first prove, in Lemma~1 of subsection 2.2, that
$\E_k'$ has no zeros of modulus $1$ in $\D$
when $k \not\equiv 2 \bmod{6}$, and that its only zeros 
of modulus $1$ in $\D$ are simple zeros at $e^{\pi i/3}$ and $e^{2\pi i/3}$
when $k \equiv 2 \bmod{6}$. This among other things implies 
assertion~$a)$ of Theorem 2. We shall then prove, in Lemma 2
of subsection 2.3, by a suitable 
application of the residue theorem, that the sum of the
multiplicities of the zeros of $\E_k'$ with modulus 
$>1$ and with real part in $]-\frac{1}{2}, \frac{1}{2}]$
is equal to $[\frac{k-4}{6}]$. The restriction of $\E_k'$ to the
half-line $\frac{1}{2} + i ]0, + \infty [$ is purely imaginary.
We shall prove in Corollary of Proposition 1
of subsection 3.1 that, with the notations of Theorem 3,
$i\E_k'$  takes non-zero values of opposite signs 
at $\frac{1}{2} + t_m$ and $\frac{1}{2} + it_{m+1}$ 
for $1 \le m \le \M-1$, and also at $\frac{1}{2} + \frac{\sqrt{3}}{2}$
and $\frac{1}{2} + it_{\M}$ when $k \equiv 4\bmod{6}$
and $k \ne 4$.  This yields us a total of $[\frac{k-4}{6}]$ disjoint open intervals
on the half-line $\frac{1}{2} + i ]\frac{\sqrt{3}}{2}, + \infty[$
in each of which $\E_k'$ has at least one zero.
But then necessarily each of these intervals can contain only 
one zero of $\E_k'$,  this zero must be simple, and there cannot 
be any other zero of $\E_k'$ whose modulus is $>1$ and whose real
part belongs to $]- \frac{1}{2}, \frac{1}{2}]$. This proves simultaneously 
Theorem 1, part $b)$ of Theorem 2 and Theorem~3.

\medskip
\subsection{Zeros of $\E_k'$ in $\gamma\D$, 
for $\gamma \in \SL_2(\Z)$ \nopunct}$~~$\\

\noindent{\bf Theorem 4.}$-$
{\em  All the zeros of $\E_k'$ in $\mathfrak{H}$ are simple.}

\medskip
We shall prove this theorem in section 4. It is a consequence of much
more general results (see subsection 4.4, Theorem 6 and 
subsection 4.5, Theorem 7) on the multiplicity of zeros of 
quasi-modular forms for $\SL_2(\Z)$ with algebraic Fourier 
coefficients, these results being themselves deduced from a 
theorem of algebraic independence of G. V. Chudnovsky.

\smallskip
Let us call {\it trivial} zeros of $\E_k'$ those zeros of $\E_k'$ 
which are also zeros of $\E_k$.  By subsection 1.2, there exist such 
zeros only when $k \equiv 2 \bmod 6$, and these are then
the elements of the orbit of $e^{\pi i/3}$
under the action of $\SL_2(\Z)$. Let us denote by $\bZ(\E_k')$ 
the set of non-trivial zeros of $\E_k'$. It is translation invariant
by the elements of $\Z$ since $\E_k'$ is periodic with period $1$.

\smallskip
We are interested in the description of the points of $\bZ(\E_k')$
belonging to $\gamma\D$, when $\gamma =\big({ a\ b \atop c\ d }\big)$
is an element of $\SL_2(\Z)$. When $c=0$, $\gamma$ is of the
form $ \pm \big({1\ r \atop 0\ 1}\big)$, $\gamma\D$ is the translate
of $\D$ by $r$, and the description of $\bZ(\E_k') \cap \gamma\D$
can be deduced from the results in subsection 1.3. For $c \ne 0$,
we have the following results : 

\bigskip
\noindent
{\bf Theorem 5.}$-$
{\em 
a) If $|d| < |c|$, the set $\bZ(\E_k') \cap \gamma\D$ is empty.

\smallskip
b) If $|d| \ge |c| > 0$, the set $\bZ(\E_k') \cap \gamma\D$ 
has cardinality $[\frac{k+2}{6}]$. It is contained in the interior of
$\gamma\D$ when $|d| \ne |c|$. When $|d| = |c|$, it is contained in
the boundary of $\gamma\D$, and more precisely in $\gamma \rC$, where
$\rC$ is the open arc of the unit circle consisting of the 
points $e^{i\theta}$, where $\frac{\pi}{3} < \theta < \frac{2\pi}{3}$. }

\medskip
Note that when $|d|=|c|$, $\gamma$ is of the 
form $\pm \big({ r+1\ r \atop \  \ 1 \ \ 1 }\big)$
or $\pm \big({\ r \  -r-1 \atop 1 \ \  -1 }\big)$
with $r \in \Z$, and $\gamma \rC$ is then the
open vertical segment 
$r + \frac{1}{2} + i ]\frac{\sqrt{3}}{6}, ~\frac{\sqrt{3}}{2}[$.

\bigskip
Theorem 5 will be proved in subsection 5.6. Here is a
consequence :

\bigskip
\noindent
{\bf Corollary.}$-$
{\em The number of zeros of $\E_k'$ in $\mathfrak{H}$ with real
part $\frac{1}{2}$ is $1 + 2[\frac{k-2}{6}]$. More
precisely, $[\frac{k-4}{6}]$ of them have an imaginary part
strictly greater than $\frac{\sqrt{3}}{2}$, $[\frac{k+2}{6}]$ 
of them have an imaginary part strictly between 
$\frac{\sqrt{3}}{6}$ and $\frac{\sqrt{3}}{2}$, none have 
an imaginary part strictly less than $\frac{\sqrt{3}}{6}$,
and when $k \equiv 2 \bmod{6}$, one has an
imaginary part equal to $\frac{\sqrt{3}}{6}$ 
and another one an imaginary part equal to $\frac{\sqrt{3}}{2}$.}

\medskip
The number of zeros of $\E_k'$ in the set
$\frac{1}{2} + i ]\frac{\sqrt{3}}{2}, + \infty[$ is 
$[\frac{k-4}{6}]$ by Theorem 2,~$b)$.

\smallskip
The set $\frac{1}{2} + i ]\frac{\sqrt{3}}{6}, \frac{\sqrt{3}}{2}[$
is equal to $\big({ 1 \ 0 \atop 1 \ 1 } \big)\rC$, with the 
notations of Theorem~5. By Theorem 5 $b)$, it 
contains $[\frac{k+2}{6}]$ zeros of
$\E_k'$, and its endpoints $\frac{1}{2} + i \frac{\sqrt{3}}{6}$
and $\frac{1}{2} + i \frac{\sqrt{3}}{2}$ are zeros of $\E_k'$
if and only if  $k \equiv 2 \bmod{6}$.

\smallskip
Finally, the set $\frac{1}{2} + i ]0, \frac{\sqrt{3}}{6}[$ is the 
image by the homography $z \mapsto \frac{z}{2z+1}$
of the half line $-\frac{1}{2} + i ]\frac{\sqrt{3}}{2}, + \infty[$. 
Then by Theorem 5 $a)$,  it does not contain any zero of $\E_k'$.

\smallskip
Thus the total number of zeros of $\E_k'$ with real part $\frac{1}{2}$
is $[\frac{k-4}{6}] + [\frac{k+2}{6}] = \frac{k -3}{3}$ when
$k \equiv 0 \bmod{6}$, $[\frac{k-4}{6}] + [\frac{k+2}{6}] + 2 = \frac{k+1}{3}$ when
$k \equiv 2 \bmod{6}$, $[\frac{k-4}{6}] + [\frac{k+2}{6}] = \frac{k-1}{3}$ when
$k \equiv 4 \bmod{6}$, hence is $1 + 2[\frac{k-2}{6}]$ in 
each of the three cases.

\bigskip
The corollary of Theorem 5 gives complete information on the number
and location of zeros of $\E_k'$ with real part $\frac{1}{2}$. Partial 
results about these zeros had previously been obtained by R. Balasubramanian 
and S. Gun (see \cite{BG}, th. 1.2 and th. 1.5).

\bigskip
\section{Counting the zeros of $\E_k'$  \\
in the standard fundamental domain}

\smallskip
As before, $k$ is an even integer $\ge 4$ and $\E_k$ is the normalized
Eisenstein series of weight $k$ for $\SL_2(\Z)$.

\medskip
\subsection{Variation of the argument of a function \nopunct}$~~$\\

Let $g$ be a continuous function on a closed bounded interval
$[a, b]$~of $\R$, with values in $\C^*$. Since $[a,b]$ is simply
connected, there exists a continuous function $h$ from
$[a, b]$ to $\R$ such that $h(t)$ is an argument of
$g(t)$ for every $t\in [a,b]$. Such a function 
is called {\it a continuous argument
of $g$.} It is unique up to addition of a constant 
function on $[a,b]$ with values in $2\pi \Z$. 
The real number $h(b) - h(a)$ does not depend on the choice 
of $h$, and is called {\it the variation of the argument 
of $g$ along~$[a, b]$}.

\bigskip
\noindent
{\bf Remark 3.}$~-$ When $g$ is of class $\rC^1$, the variation of its argument 
along $[a, b]$ is the imaginary part of 
$\int_{a}^{b} \frac{g'(t)}{g(t)}~dt$.

\medskip
\noindent
{\bf Remark 4.}$~-$ Let $c$ be a point of $[a, b]$. The variation
of the argument of $g$ along $[a, b]$ is the sum
of the variations of the argument of $g$ along 
$[a, c]$ and along $[c, b]$.

\medskip
\noindent
{\bf Remark 5.}$~-$ Let $[\alpha,\beta]$ be another
closed bounded interval of $\R$ and let $u : [\alpha, \beta] \to[a,b]$
be a continuous function which maps $\alpha$ to $a$ and
$\beta$ to $b$. The variation of the argument of 
$g$ along $[a,b]$ is equal to that of $g \circ u$ along~$[\alpha, \beta]$.

\medskip

More generally, let $\gamma : [a,b] \to \C$ be a continuous map 
and let $\varphi$ be a continuous function on $\gamma([a,b])$,
with values in $\C^*$. We define {\it variation of the argument of }
$\varphi$ {\it along the oriented curve} $\gamma$ to be the variation
of the argument of $\varphi \circ \gamma$ along $[a,b]$. By Remark 5,
it is invariant by a continuous change of parametrization of the oriented curve.

\medskip 
\subsection{Zeros of modulus $1$ of $\E_k'$ in $\D$ \nopunct}$~~$\\

\noindent{\bf Lemma 1.}$-$
{\em When $k \not\equiv 2 \bmod 6$, $\E_k'$ has no
zeros of modulus $1$ in $\D$. When
$k \equiv 2 \bmod 6$, its only zeros of modulus $1$
in $\D$ are simple zeros at $e^{\pi i/3}$ and $e^{2\pi i/3}$.}

\smallskip
For all $z \in \gH$, we have
$\E_k(- \overline{z}) = \overline{\E_k(z)}$, 
and hence 
\begin{equation}\label{eq1}
\E_k'(-\overline{z}) = - \overline{ \E_k' (z) }.
\end{equation}
Since $\E_k$ is a modular form of weight $k$ 
for $\SL_2(\Z)$, we also have
$\E_k(-\frac{1}{z}) = z^k \E_k(z),$
and hence 
\begin{equation}\label{feq}
\E_k'(-\frac{1}{z}) = z^{k+2}\E_k'(z) + k z^{k+1}\E_k(z).
\end{equation}
If $z$ is a zero of $\E_k'$ of modulus $1$, 
$-\overline{z}$ is equal to $-\frac{1}{z}$, and
we then deduce from \eqref{eq1} 
and \eqref{feq} that $\E_k(z) = 0$.
Hence $z$ is a multiple zero of $\E_k$.
But, by subsection 1.2, $\E_k$ does not have 
multiple zeros in $\D$ unless $k \equiv 2 \bmod 6$,
and these are then zeros of order $2$ of $\E_k$,
hence simple zeros of $\E_k'$, located  
at $e^{\pi i/3}$ and $e^{2\pi i/3}$.

\medskip
\subsection{Other zeros of $\E_k'$ in $\D$ \nopunct}$~~$\\

We have determined in subsection 2.2 the zeros of 
$\E_k'$ in $\D$ which are of modulus $1$. 
Let us now look at the other zeros of $\E_k'$ in $\D$, namely
those of modulus $> 1$. Since those with real part $\frac{1}{2}$
are translates by $1$ of those with real part $-\frac{1}{2}$, 
we will only be interested in those zeros whose real 
part is distinct from $-\frac{1}{2}$.
The following lemma allows us to count them :

\bigskip
\noindent{\bf Lemma 2.}$-$
{\em The sum of the multiplicities of those zeros 
of $\E_k'$ with modulus $>1$ and real part in 
$] -\frac{1}{2}, \frac{1}{2}]$ is equal to $[\frac{k-4}{6}]$.}

\medskip
When $\Im(z)$ tends to $+ \infty$, 
$\E_k'(z)$ is equivalent to $-\frac{4 k\pi i }{\B_k} e^{2 \pi i z}$.
Therefore we can choose a real number $\T > 1$ such that
all zeros of $\E_k'$ have an imaginary part 
strictly less than $\T$. By the residue theorem, the sum
considered in Lemma~2 is equal, for $\varepsilon > 0$ sufficiently 
small, to the integral
\begin{equation}\label{cz}
\I :=  \frac{1}{2\pi i}\int_{\gamma_{\T, \varepsilon}}
\frac{\E_k''(z)}{\E_k'(z)} ~dz,
\end{equation}
where $\gamma_{\T, \varepsilon}$ is the contour of integration
obtained as follows : one modifies the 
oriented boundary of $\D_{\T} = \{ z \in \D ~|~ \Im(z) \le {\T}\}$
by going around the zeros of $\E_k'$ lying on the boundary
of $\D_{\T}$ by arcs of circles of radius $\varepsilon$
passing outside of $\D_{\T}$ for the zeros with modulus $>1$ 
and real part $\frac{1}{2}$, and passing inside  $\D_{\T}$ 
for all the others (see Figure 1 below \footnote[1]{\,\tiny{The authors thank Sunil Naik, 
a doctoral student at the Institute of Mathematical Sciences in Chennai, who 
carried out the drawing of all figures in this paper.}\par}).

\bigskip

\centerline{\includegraphics[width=9.6cm, height=8.4cm]{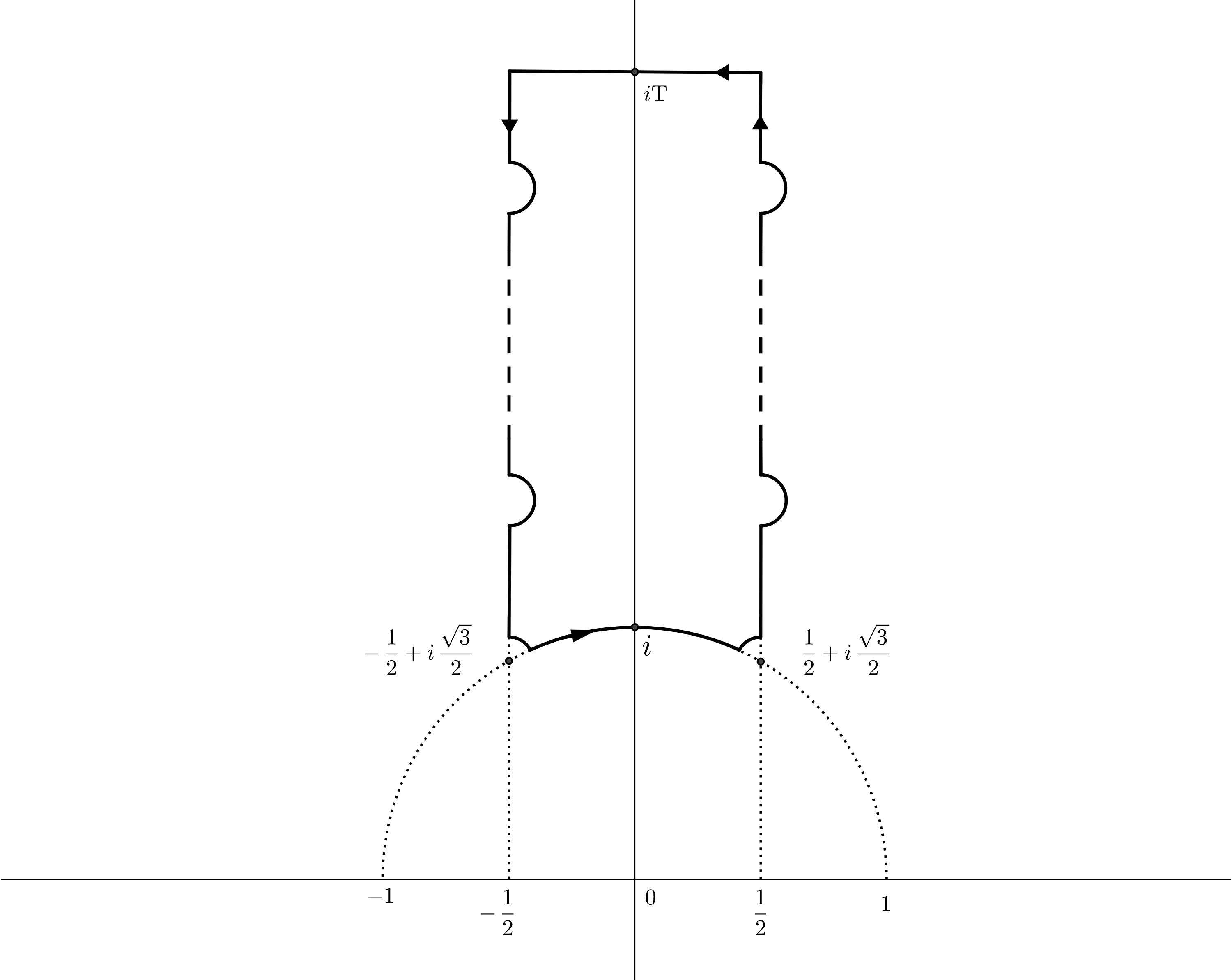}}

\medskip
\centerline{{\it Figure 1.} 
The contour of integration  $\gamma_{\T,\varepsilon}$ 
when $k\equiv 2\bmod 6$.}

\bigskip
The computation of the integral $\I$ will be the subject of the
remaining parts of this section. We will show that $\I = [\frac{k-4}{6}]$
by distinguishing three cases according to the 
congruence class of $k$ modulo $6$ (see Lemmas 3, 4
and 5). Lemma 2 then will follow.

\medskip
\subsection{Reduction of the computation of the integral $\I$ to 
a variation of argument \nopunct}$~~$\\

The contribution to the integral $\I$ of the upper horizontal 
oriented segment (joining $\frac{1}{2} + i \T$ to $-\frac{1}{2} + i \T$) 
is independent of $\T$ since 
$\E_k$ is periodic with period~$1$. This contribution is equal 
to $-1$ : this is seen by letting $\T$ go to
$+\infty$, and observing that $\frac{\E_k''(z)}{\E_k'(z)}$
tends to $2\pi i$ when $\Im(z)$ goes to~$+ \infty$.

\smallskip
The contributions to the integral $\I$ of the vertical segments
contained in the line with abscissa $\frac{1}{2}$ and 
of those contained in the line with abscissa $-\frac{1}{2}$
are opposite of each other, since $\E_k$ is periodic with period $1$
and these segments have opposite orientations.

\smallskip
Similarly, the contributions to the integral $\I$
of the semicircles of radius $\varepsilon$ centered at
the zeros of $\E_k'$ of modulus $>1$ and real part
$\frac{1}{2}$ are opposite to those of the semicircles
of radius $\varepsilon$ centered at
the zeros of $\E_k'$ of modulus $>1$ and real part
$-\frac{1}{2}$.

\smallskip
When $k \not\equiv 2 \bmod{6}$, $\E_k'$ has no 
zeros of modulus $1$ in $\D$ by Lemma~1.
The only remaining contribution to the integral $\I$
is then that of the arc of the unit circle joining 
$e^{2\pi i/3}$ to $e^{\pi i/3}$.
It can be written as $-\frac{\A}{2\pi}$, where
\begin{equation}\label{eq2.1}
\A
:= 
\int_{\frac{\pi}{3}}^{\frac{2\pi}{3}} 
\frac{\E_k''(e^{i\theta})}{\E_k'(e^{i\theta})} 
e^{i \theta} ~d\theta.
\end{equation}
We have $\I = -1 - \frac{\A}{2\pi}$.
Since $\I$ is an integer, $\A$ is a real number. 
It follows that $\A$ is the variation of the argument 
of $\theta \mapsto \E_k'(e^{i \theta})$ 
along the interval $\left[\frac{\pi}{3}, \frac{2\pi}{3}\right]$
(see 2.1, Remark 3).

\smallskip
When $k \equiv 2 \bmod 6$, it follows from Lemma 1 that the 
only zeros of $\E_k'$ in $\D$ of modulus~$1$ are simple zeros
at $e^{\pi i/3}$ and $e^{2\pi i/3}$. The function
$\frac{\E_k''(z)}{\E_k(z)}$ is equivalent to 
$\frac{1}{z- e^{\pi i/3}}$ when $z$ approaches
$e^{\pi i/3}$. The contribution to $\I$
of the circular arc of radius $\varepsilon$ centered
at $e^{\pi i/3}$ tends to $-\frac{1}{6}$ when $\varepsilon$
tends to $0$. The same holds for the contribution
to $\I$ of the circular arc of radius $\varepsilon$
centered at $e^{2 \pi i/3}$.
Hence we have $\I =  - \frac{4}{3} - \frac{\A}{2\pi}$,
where 
\begin{equation}\label{lvar}
\A
:= 
\lim_{\eta \to 0^+} 
\int_{\frac{\pi}{3} + \eta}^{\frac{2\pi}{3} - \eta}
\frac{\E_k''(e^{i\theta})}{\E'_k(e^{i\theta})}
e^{i\theta}~d\theta.
\end{equation}
We see as before that $\A$ is real and is equal to the limit, when
$\eta >0$ tends to $0$, of the variation of the 
argument of 
$\theta \mapsto \E_k'(e^{i\theta})$
along the interval 
$\left[ \frac{\pi}{3} + \eta,  \frac{2\pi}{3} - \eta \right]$.

\medskip
\subsection{Behaviour of $\E_k'$ on the unit circle \nopunct}$~~$\\

For any point $z$ in $\gH$ of modulus $1$, we have
\begin{equation}\label{2.5.1}
\overline{\E_k(z)} ~=~ \E_k(-\bar{z}) 
~=~ \E_k(-\frac{1}{z}) ~=~ z^k \E_k(z).
\end{equation}
Hence, for all $\theta \in ]0, \pi[$, we can write
\begin{equation}\label{eqN2.2}
\E_k(e^{i\theta}) = e^{-k i \theta/2} f_k(\theta),
\end{equation}
where $f_k(\theta)$ is a real number.  The function $f_k$ defined 
in this manner on $]0, \pi [$ is real analytic. By differentiating, we get
\begin{equation}\label{eq2.2}
\E_k'(e^{i\theta}) 
~=~ 
-i e^{-(k+2) i \theta/2}g_k(\theta),
\end{equation}
where
\begin{equation}\label{neqfn}
g_k(\theta) := f_k'(\theta) - \frac{k i}{2} f_k(\theta).
\end{equation}

\medskip
\noindent
{\bf Remark 6.}$~-$ It follows from formulas
\eqref{2.5.1} and \eqref{eqN2.2} that
we have $f_k(\pi -\theta) = (-1)^{k/2} f_k(\theta)$
for $\theta \in ]0, \pi [$. So we have
$f_k'(\pi -\theta) = -(-1)^{k/2}f_k'(\theta)$, and
$g_k(\pi - \theta)$ is symmetric of $g_k(\theta)$
 with respect to the horizontal axis when
 $k \equiv 2\bmod{4}$, and to the
 vertical axis when $k \equiv 0\bmod{4}$.

\medskip
It follows from Lemma 1 of subsection 2.2 that $g_k$ has no 
zeros in the interval $[\frac{\pi}{3}, \frac{2\pi}{3}]$
when $k \not\equiv 2 \bmod{6}$, and has 
zeros only at the endpoints of this interval
when $k \equiv 2\bmod{6}$. We then deduce from 
\eqref{eq2.2} that
\begin{equation}\label{neweq2.5}
\A 
~=~ 
-\frac{k+2}{6}\pi  + \B,
\end{equation}
where $\B$ is the variation of the argument of 
$g_k$ along $[\frac{\pi}{3}, \frac{2\pi}{3}]$
when $k \not\equiv 2 \bmod{6}$, and the limit,
when $\eta >0$ tends to $0$, of the 
variation of the argument of $g_k$ along
$[\frac{\pi}{3}+ \eta,  \frac{2\pi}{3} -\eta]$
when $k \equiv 2\bmod{6}$.

\medskip
\subsection{Computation of the integral $\I$ when 
$k \equiv 4 \bmod 6$ \nopunct}$~~$\\

\noindent{\bf Lemma 3.}$-$
{\em When $k \equiv 4\bmod 6$, the integral $\I$ is equal to $\frac{k-4}{6}$.}

\medskip
Suppose that $k \equiv 4 \bmod{6}$. It follows from subsection 1.2
that the function $f_k$ has only simple
zeros in the interval $[\frac{\pi}{3}, \frac{2\pi}{3}]$, and that they
can be written as 
\begin{equation}\label{eqN2.6}
\theta_0 < \theta_1 < \cdots < \theta_{\rN} < \theta_{\rN+1}
\end{equation}
where $\theta_0 = \frac{\pi}{3},~\theta_{\rN+1} = \frac{2\pi}{3}$,
and $\rN$ is equal to $2[\frac{k}{12}]$
when $k \equiv 4 \bmod{12}$ and to
$2[\frac{k}{12} ] +1$
when $k \equiv 10 \bmod{12}$, and hence
to $\frac{k-4}{6}$ in both of these cases.

\medskip
The points $\theta_j$ are those points in the interval
$[\frac{\pi}{3}, \frac{2\pi}{3}]$,
where the complex-valued function $g_k = f_k' - \frac{ki}{2}f_k$
crosses the real line. Since $\theta_j$ is a simple zero
of $f_k$, the real number $g_k(\theta_j) = f_k'(\theta_j)$
is non-zero. If it is positive, $g_k$ crosses the real line at $\theta_j$ 
from upper half-plane to lower half-plane. If it is negative,
$g_k$ crosses the real line at $\theta_j$ from lower half-plane to 
upper half-plane. Therefore $g_k(\theta_j)$ and $g_k(\theta_{j+1})$ 
have opposite signs for $0 \le j \le \rN$, and more precisely :

\smallskip
$a)$ If $g_k(\theta_j) < 0$, $g_k([\theta_j, \theta_{j+1}])$ is contained
in the closed upper half-plane, and $g_k$ has a 
continuous argument along $[\theta_j, \theta_{j+1}]$
taking the values $\pi$ at $\theta_j$ and $0$ at~$\theta_{j+1}$;

\smallskip
$b)$ If $g_k(\theta_j) > 0$, $g_k([\theta_j, \theta_{j+1}])$ is contained
in the closed lower half-plane, and $g_k$ has a 
continuous argument along $[\theta_j, \theta_{j+1}]$
taking the values $0$ at $\theta_j$ and $-\pi$ at~$\theta_{j+1}$.

\smallskip
We observe that the variation of the argument of $g_k$
along each interval $[\theta_j, \theta_{j+1}]$, where  $0 \le j \le \rN$,
is equal to $-\pi$. Hence the variation $\B$ of the argument of $g_k$
along the interval $[\frac{\pi}{3}, \frac{2\pi}{3}]$ is equal to
$-(\rN +1)\pi = - \frac{k+2}{6}\pi$ (subsection 2.1, Remark 4).
It follows from formula \eqref{neweq2.5} that 
$\A = -\frac{k+2}{6} \pi~+~\B~=~-~\frac{k+2}{3}\pi$, and then 
by subsection 2.4 that
$\I = -1 - \frac{\A}{2\pi} = \frac{k-4}{6}$.

\bigskip
\noindent
{\bf Remark 7.}$~-$ Although it is not required for the paper, 
let us point out that for $0 \le j \le \rN+1$, the sign of 
$g_k(\theta_j)$ is equal to $(-1)^{ \rN + 1 -j}$. Since we already know
that consecutive terms of this sequence have opposite signs, it is 
sufficient to show that $g_k(\theta_{\rN +1}) = g_k(\frac{2\pi}{3})$ is positive.
Now, by subsection 1.2, the real analytic function
$t \mapsto \E_k(-\frac{1}{2} + it)$
has a simple zero at $\frac{\sqrt{3}}{2}$ and is strictly positive 
for $t > \frac{\sqrt{3}}{2}$. It follows that its derivative at the point 
$\frac{\sqrt{3}}{2}$ is strictly positive.  This derivative is
$i \E_k'(e^{2\pi i/3})$. It is equal to $g_k(\frac{2\pi}{3})$
by formula~\eqref{eq2.2}, since $k \equiv 4 \bmod 6$.

\bigskip

\centerline{\includegraphics[width=9.6cm, height=6cm]{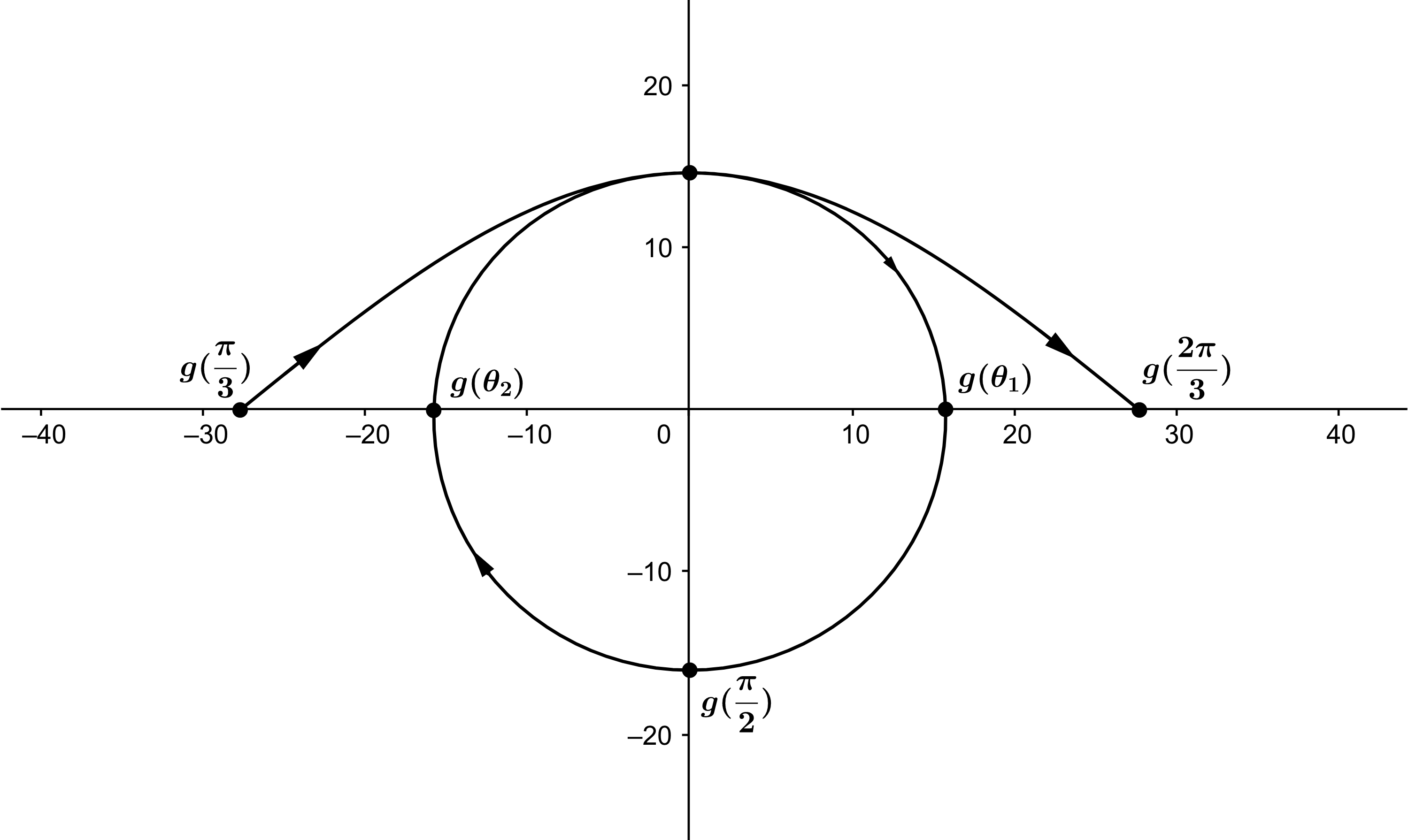}}

\medskip
\centerline{{\it Figure 2.}  The oriented curve $g([{\pi\over3},{2\pi\over 3}])$, 
where $g=g_{16}$.}

\medskip
\subsection{Computation of the integral $\I$ when $k \equiv 0 \bmod 6$ \nopunct}$~~$\\

\noindent{\bf Lemma 4.}$-$
{\em When $k \equiv 0 \bmod 6$, the integral $\I$ is equal to 
$\frac{k-6}{6}$, and so to $[\frac{k-4}{6}]$.}

\medskip
Suppose that $k \equiv 0 \bmod{6}$. It follows from subsection 1.2
that the zeros of $f_k$ in the interval $[\frac{\pi}{3}, \frac{2\pi}{3}]$ are simple, 
distinct from $\frac{\pi}{3}$ and $\frac{2\pi}{3}$, and that their number 
$\rN$ is equal to $2[\frac{k}{12}]$
when $k \equiv 0 \bmod{12}$ and to $2[\frac{k}{12}] +1$
when $k \equiv 6 \bmod{12}$, i.e. to $\frac{k}{6}$ in both
cases. Let us write these zeros as
\begin{equation}\label{eqN2.7}
\theta_1 < \cdots < \theta_{\rN}.
\end{equation}

\smallskip
One proves as in subsection 2.6 that the $g_\k(\theta_j)$ are real,
non-zero, of alternating signs,
and that the variation of the argument of $g_\k$ along 
$[\theta_j, \theta_{j+1}]$ is $-\pi$ for $1 \le j \le \rN-1$.

\smallskip
The integer $\rN = \frac{k}{6}$ is non-zero. If 
$g_k(\theta_{\rN}) >0$ (resp. $g_k(\theta_{\rN}) < 0$), the
set $g_k(]\theta_{\rN}, \frac{2\pi}{3}])$ is contained in the open
lower (resp. upper) half-plane. 
In both of these cases, the
variation of the argument of $g_k$ along the interval 
$[\theta_{\rN}, \frac{2\pi}{3}]$ is strictly between $-\pi$ and $0$.
We show in a similar way (or we deduce from Remark 6)
that the variation of the argument
of $g_k$ along the interval $[\frac{\pi}{3}, \theta_1]$
is also strictly between $-\pi$ and $0$.

\smallskip
Summing up these contributions, we
see that the variation $\B$ of the argument
of $g_k$ along the interval $[\frac{\pi}{3}, \frac{2\pi}{3}]$
is strictly between $-(\rN+1)\pi$ and $-(\rN -1)\pi$,
i.e. between $- \frac{k+6}{6} \pi$ and $-\frac{k-6}{6}\pi$.
Since $\A = - \frac{k+2}{6}\pi + \B$ by formula \eqref{neweq2.5}
and $\I = -1 - \frac{\A}{2\pi}$ by subsection 2.4, $\I$  is strictly between 
$\frac{k-8}{6}$ and $\frac{k-2}{6}$.
We know that $\I$ is an integer and by hypothesis we have
$k \equiv 0 \bmod{6}$. It follows that $\I ~=~ \frac{k-6}{6}$ (and 
consequently $\A = - \frac{k}{3}\pi$ and $\B = - \frac{k-2}{6}\pi$).

\bigskip
\noindent
{\bf Remark 8.}$-$ Although it is not required, let us list few additional 
details. By subsection 1.2, the real analytic function $t \mapsto \E_k(\frac{1}{2} + it)$
is strictly positive for $t > \frac{\sqrt{3}}{2}$
and is non-zero at $\frac{\sqrt{3}}{2}$. We hence have
$\E_k (e^{2\pi i/3}) = \E_k(e^{\pi i/3}) > 0$. Therefore, using 
formula \eqref{eqN2.2},
we have $f_k(\frac{2\pi}{3}) > 0$, the sign of $f_k(\frac{\pi}{3})$ 
is $(-1)^{k/6}$, and the sign of $g_k(\theta_j)$ is $(-1)^{\rN -j}$ for
$1 \le j \le \rN$.

\bigskip
\noindent
{\bf Remark 9.}$-$ The restriction of $\E_k'$ to the half-line 
$-\frac{1}{2} + i]0, +\infty[$ is purely imaginary. We therefore deduce from
formula \eqref{eq2.2} that $g_\k(\frac{2\pi}{3})$ 
belongs to $e^{2\pi i/3}\R$, and in fact to $e^{-\pi i/3}\R_+^*$ since its
imaginary part is strictly negative (Remark 8). It follows that the variation of
the argument of $g_k$ along the interval $[\theta_{\rN}, \frac{2\pi}{3}]$ is
$-\frac{\pi}{3}$. We show in a similar way (or deduce from Remark 6) that the
variation of the argument of $g_k$ along the interval $[\frac{\pi}{3}, \theta_1]$
is also equal to $-\frac{\pi}{3}$.  This provides a new proof, more direct 
than the previous one, of the fact that the variation of the argument of $g_k$ 
along the interval $[\frac{\pi}{3}, \frac{2\pi}{3}]$ is equal to 
$-(\rN -1)\pi - \frac{2\pi}{3}$,  i.e. to $-\frac{k-2}{6}\pi$, and so of Lemma 4.

\bigskip
 
\centerline{\includegraphics[width=6cm, height=6cm]{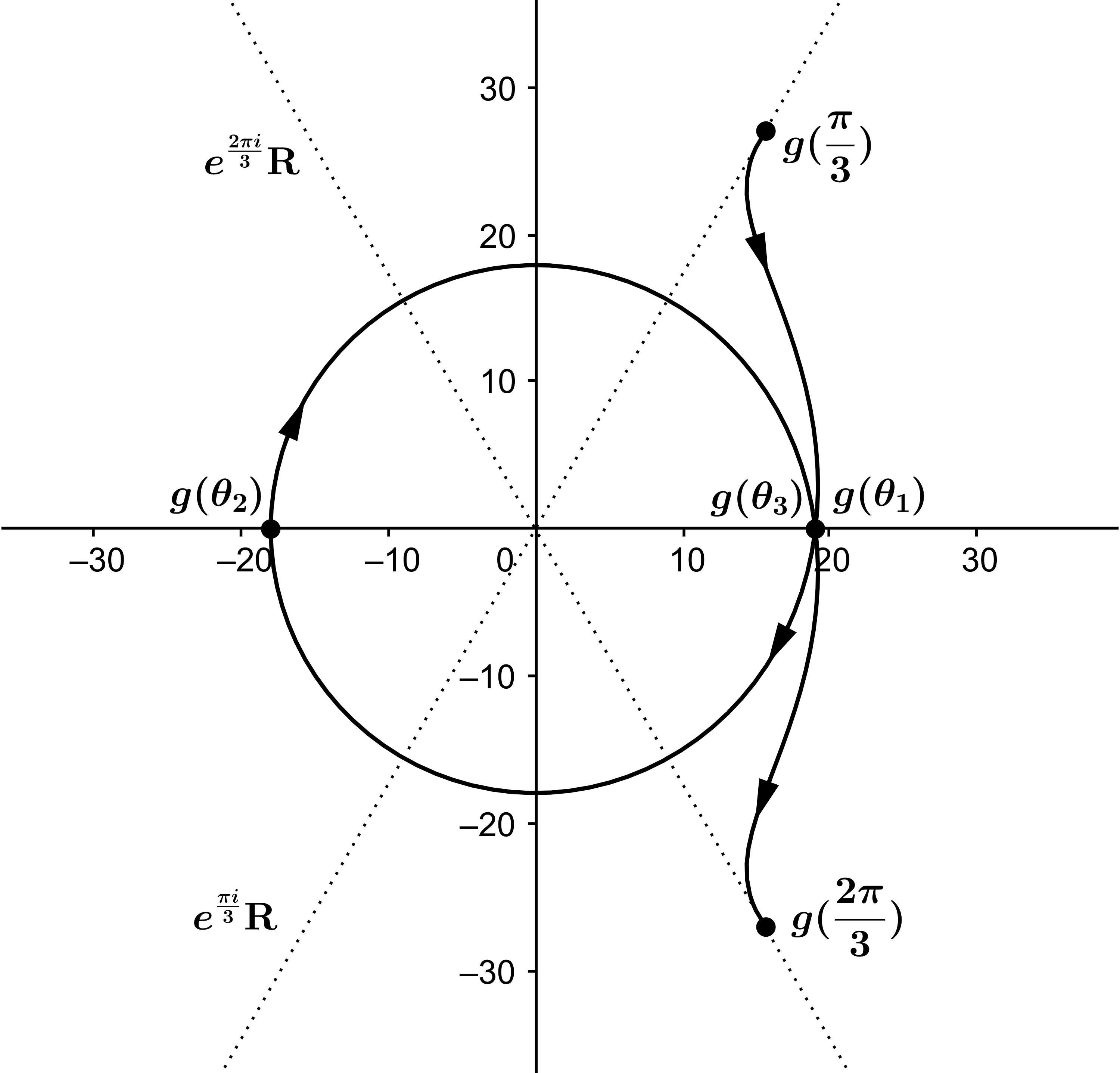}}

\medskip

\centerline{{\it Figure 3.}  The oriented curve $g([{\pi\over3},{2\pi\over 3}])$, 
where $g=g_{18}$.}\

\medskip
\subsection{Computation of the integral $\I$ when $k \equiv 2 \bmod 6$\nopunct}$~~$\\

\noindent{\bf Lemma 5.}$-$
{\em When $k \equiv 2\bmod 6$, the integral $\I$ is equal to $\frac{k-8}{6}$, and 
so to $[\frac{k-4}{6}]$.}

\medskip
Suppose that $k \equiv 2 \bmod{6}$. It follows from subsection 1.2
that $f_k$ has zeros of order $2$ at $\frac{\pi}{3}$ and $\frac{2\pi}{3}$,
that the zeros of $f_k$ in the interval $]\frac{\pi}{3}, \frac{2\pi}{3}[$ 
are simple, and that
their number $\rN$ is equal to $2[\frac{k}{12}]$
when $k \equiv 8 \bmod{12}$ and to $2[\frac{k}{12}] -1$
when $k \equiv 2 \bmod{12}$, i.e. to $\frac{k-8}{6}$ in both
the cases. Let us write the latter ones as
\begin{equation}\label{eqN2.8}
\theta_1 < \cdots < \theta_{\rN}.
\end{equation}
One proves as in subsection 2.6 that the $g_k(\theta_j)$ are real, non-zero, 
of alternating signs, and that the variation of the argument of $g_k$ 
along $[\theta_j, \theta_{j+1}]$ is $-\pi$ for $1 \le j \le \rN -1$.

\smallskip
When $\rN \ne 0$ (i.e. $k \ne 8$), let us distinguish two cases :

\smallskip
$a)$ If $g_\k(\theta_{\rN}) > 0$, $g_k([\theta_{\rN}, \frac{2\pi}{3}])$ is 
contained in the closed lower half-plane. Let $h$ be the 
continuous argument of $g_k = f_k' - \frac{k i}{2}f_\k$ in 
$[\theta_{\rN}, \frac{2\pi}{3}[$ which takes the 
value $0$ at $\theta_{\rN}$. When $\theta < \frac{2\pi}{3}$
tends to $\frac{2\pi}{3}$, $f_\k(\theta)$ is equivalent to 
$c(\theta - \frac{2\pi}{3})^2$ for some positive real number 
$c$, $f_\k'(\theta)$ is equivalent to 
$2c(\theta - \frac{2\pi}{3})$, hence the real and imaginary parts of 
$g_\k(\theta)$ are strictly negative and their quotient 
$\frac{\Im \,g_\k(\theta)}{\Re\, g_\k(\theta)}$ tends to $0$, which implies 
that $h(\theta)$ tends to $-\pi$.

\smallskip
$b)$ If $g_k(\theta_{\rN}) < 0$, $g_k([\theta_{\rN}, \frac{2\pi}{3}])$ is 
contained in the closed upper half-plane. Let $h$ be the 
continuous argument of $g_k = f_k' - \frac{k i}{2}f_k$ in 
$[\theta_{\rN}, \frac{2\pi}{3}[$ which takes the 
value $\pi$ at $\theta_{\rN}$. When $\theta < \frac{2\pi}{3}$
tends to $\frac{2\pi}{3}$, $f_k(\theta)$ is equivalent to 
$c(\theta - \frac{2\pi}{3})^2$ for some negative
real number $c$, $f_k'(\theta)$ is equivalent to 
$2c(\theta - \frac{2\pi}{3})$, hence the real and imaginary parts of 
$g_k(\theta)$ are strictly positive and their quotient 
$\frac{\Im \,g_k(\theta)}{\Re \,g_k(\theta)}$ tends to $0$, which implies 
that $h(\theta)$ tends to~$0$.

We observe that in both the previous cases, the limit 
when $\eta >0$ tends to~$0$ of the 
variation of the argument of $g_k$ along the 
interval $[\theta_{\rN}, \frac{2\pi}{3} - \eta]$ is equal to $-\pi$. 
Similarly we prove (or we deduce from Remark 6) that
the limit when $\eta >0$ tends to $0$
of the variation of the argument of $g_k$ along the interval 
$[\frac{\pi}{3}+ \eta, \theta_1]$ is equal to $-\pi$.

Finally, when $\rN =0$ (i.e. $k=8$), 
$g_k([\frac{\pi}{3}, \frac{2\pi}{3}])$
is either contained in the closed lower half-plane
or in the closed upper half-plane. In the first (resp. the
second) case, we check as above that there is a
continuous argument $h$ of $g_k$ along the 
interval $]\frac{\pi}{3}, \frac{2\pi}{3}[$ which
has limit $0$ at $\frac{\pi}{3}$ and $-\pi$
at $\frac{2\pi}{3}$ (resp. $\pi$ at $\frac{\pi}{3}$
and $0$ at $\frac{2\pi}{3}$). The limit when
$\eta>0$ tends to $0$ of the variation of argument
of $g_k$ along the interval 
$[\frac{\pi}{3} + \eta, \frac{2\pi}{3} - \eta]$ is therefore 
equal to $-\pi$.

Summing up all these contributions, we see that, in all
cases, the limit $\B$ when $\eta >0$ tends to $0$
of the variation of the argument of $g_k$ along 
the interval $[\frac{\pi}{3} + \eta, \frac{2\pi}{3} - \eta]$
is equal to $-(\rN+1)\pi$, i.e. to $-\frac{k-2}{6}\pi$.
Since $\A = -\frac{k+2}{6}\pi + \B$ by formula
\eqref{neweq2.5} and $\I = - \frac{4}{3} - \frac{\A}{2\pi}$
by subsection 2.4, we have $\A = - \frac{k}{3}\pi$
and $\I = \frac{k-8}{6}$.

\medskip
\noindent
{\bf Remark 10.}$-$ Although it is not required for the paper, let us 
add few more clarifications. When $\k \equiv 2 \bmod{6}$, $\E_\k$ 
has by subsection 1.2 a zero of order $2$ at 
$e^{2\pi i/3}~=~-\frac{1}{2}~+~i\frac{\sqrt{3}}{2}$. Hence there exists
$c \in \C^*$ such that $\E_\k(z)$ is equivalent to $-c(z- e^{2\pi i/3})^2$ when
$z$ tends to $e^{2\pi i/3}$. Since $\E_\k(- \frac{1}{2} + it) >0$ 
for $t > \frac{\sqrt{3}}{2}$ (see subsection 1.2, Remark 2), $c$ is a positive
real number. But then $\E_\k(e^{i\theta})$ is equivalent to 
$ce^{4\pi i/3}(\theta - \frac{2\pi}{3})^2$ when $\theta$ tends to 
$\frac{2\pi}{3}$, and $f_\k(\theta)$ to $c(\theta - \frac{2\pi}{3})^2$
by formula \eqref{eqN2.2}. We then deduce from the above discussion
that $g_\k(\theta_{\rN}) >0$ (if $\rN \ne 0$), so that $g_\k(\theta_i)$
has sign $(-1)^{\rN -i}$ for $1 \le i \le \rN$.

\bigskip
 
\centerline{\includegraphics[width=7.2cm, height=7.2cm]{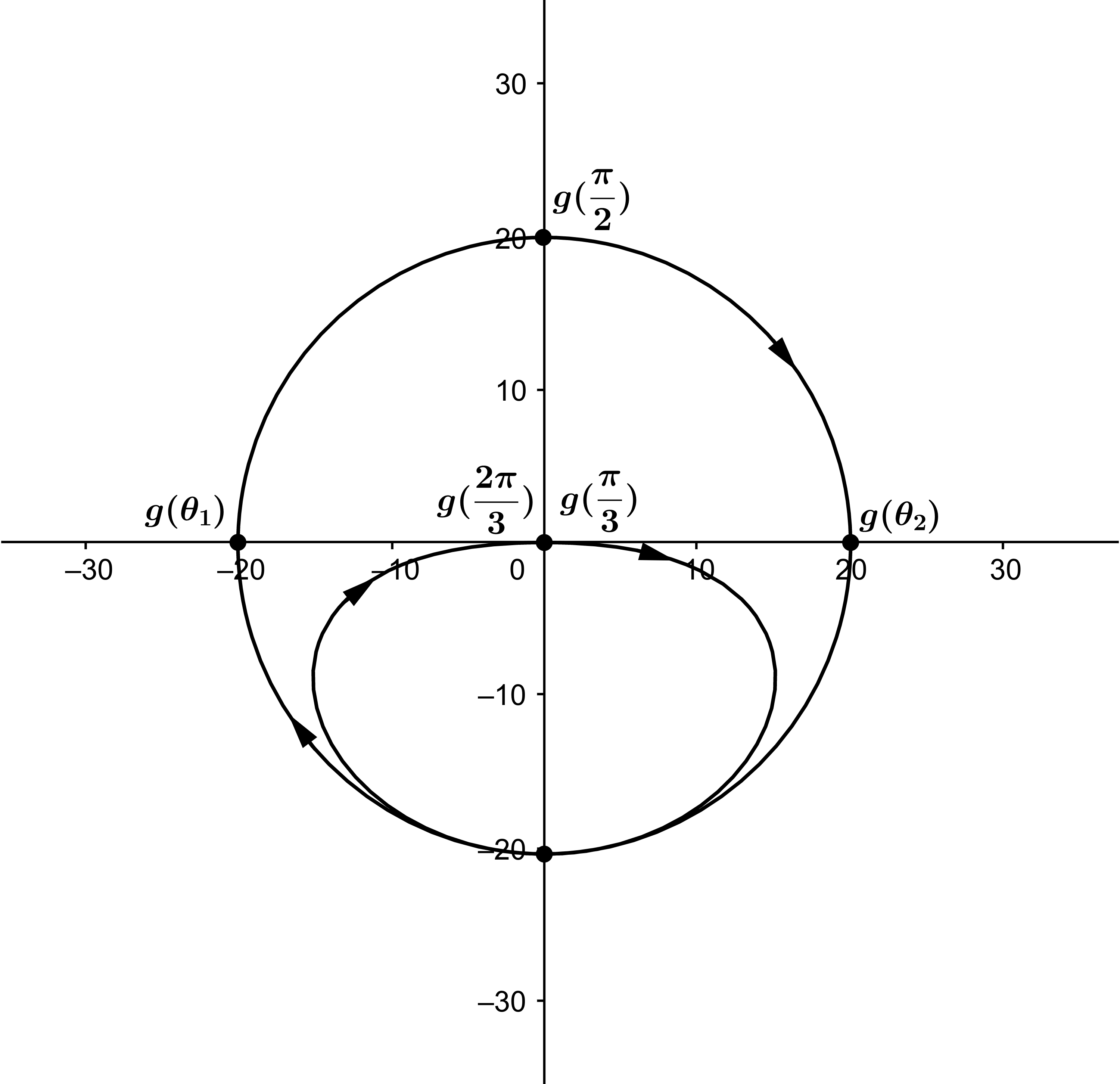}}

\medskip

\centerline{{\it Figure 4.}  The oriented curve 
$g([{\pi\over3},{2\pi\over 3}])$, where $g=g_{20}$.}

\medskip
\section{Signs changes of $\E_k'$}

\smallskip
\subsection{Statement of results \nopunct}$~~$

\smallskip
In this section, $k$ denotes an integer $\ge 2$ (not necessarily even)
and $h_k$ the function defined on the Poincar\'e upper half-plane
by
\begin{equation}\label{nf}
h_k(z) 
:= 
\sum_{n =1}^{\infty} n^k \frac{q^n}{(1- q^n)^2},
\end{equation}
where $q = e^{2\pi i z}$. Note that we have
\begin{equation}\label{vf}
h_k(z)
~=~
\sum_{n =1}^{\infty} n^k (\sum_{m =1}^{\infty} mq^{mn})
~=~
\sum_{n=1}^{\infty} n\sigma_{k-1}(n)q^n.
\end{equation}
It follows by differentiating formula \eqref{NEis}
that, when $k$ is an even integer $\ge 4$,
\begin{equation}\label{DEis}
\E_k' (z) ~=~ -\frac{4\pi k i}{\B_k} h_k(z).
\end{equation}
For $1 \le m \le \left[ \frac{k+1}{6}\right]$, let
$t_m := \frac{1}{2} \cot (\frac{m\pi}{k+1})$.
This section is devoted to prove the 
following proposition.

\medskip
\noindent{\bf Proposition 1.}$-$
{\em For any integer $m$ between $1$ and 
$\left[\frac{k+1}{6} \right]$, $h_k(\frac{1}{2} + i t_m)$
is a non-zero real number whose sign is $(-1)^m$.}

\bigskip
Note that the inequalities $1 \le m \le \left[\frac{k+1}{6} \right]$ imply that
$t_m \ge \frac{\sqrt{3}}{2}$. In subsection 3.2, we prove Proposition 1
when $3m^2 \le k$ using the expression \eqref{nf} of the function $h_k$.
We will give in subsection 3.3 a new expression of $h_k$ 
and we shall use it to prove Proposition 1 when $3m^2 \ge k+1$ 
in subsection 3.4.

\bigskip
\noindent
{\bf Corollary.}$-$
{\em 
Suppose that $k$ is an even integer $\ge 4$ and 
let $\M := [\frac{k}{6}]$. \par

$a)$ For any integer $m$ between $1$ and $\M -1$, $i\E_k'$ takes
non-zero real values of opposite signs at $\frac{1}{2} + t_m$ and 
$\frac{1}{2} + it_{m+1}$. \par

$b)$ When $k \equiv 4 \bmod 6$ and $k \ne 4$, $i\E_k'$ takes at 
$\frac{1}{2} + i \frac{\sqrt{3}}{2}$ and $\frac{1}{2} + i t_{\M}$
non-zero real values of opposite signs.
}

\smallskip
Assertion $a)$ follows from Proposition 1, using formula \eqref{DEis}.
 
\smallskip
When $k\equiv 4 \bmod{6}$ and $k \ne 4$, we have $k = 6\M + 4$ and $\M \ge 1$. 
The sign of $h_k(\frac{1}{2} + it_{\M})$ is $(-1)^{\M}$ by Proposition 1.
The sign of $\B_k$ is $-1$ or $1$ depending on whether $k$ is congruent to
$0$ or $2$ modulo $4$; it is hence equal to $(-1)^{\M -1}$. We then deduce
from formula \eqref{DEis} that $i \E_k'$ takes at $\frac{1}{2} + it_{\M}$ a strictly
negative real value. Note that the derivative at the point $\frac{\sqrt{3}}{2}$
of the function $t \to \E_k(\frac{1}{2} + it)$ is real and strictly positive
by Remark 7 of subsection 2.6. In other words, we have 
$i \E_k'(\frac{1}{2} + i \frac{\sqrt{3}}{2}) > 0$. This proves assertion $b)$.

\medskip
\subsection{Proof of Proposition 1 when $3m^2 \le k$ \nopunct}$~~$\\

Throughout this subsection, we shall assume that $m$ is an integer
such that 
$1 \le m \le \left[ \frac{k+1}{6} \right]$ and $3m^2 \le k$. We put
$t_m := \frac{1}{2} \cot (\frac{m\pi}{k+1})$, $z_m := \frac{1}{2} + it_m$
and $q_m := e^{2\pi i z_m}= -e^{-2\pi t_m}$.

\medskip
\noindent{\bf Lemma 6.}$-$ 
{\em We have
$|q_m| \le e^{-\frac{k+1}{m} + \frac{4m}{k+1}}$.} 

\medskip
For all $x \in ]0, \frac{\pi}{2}]$, we have
\begin{equation}\label{excot}
\cot x 
~=~ 
\frac{1}{x} + \sum_{n=1}^{\infty} (\frac{1}{x-n\pi} + \frac{1}{x+n\pi})
~=~ 
\frac{1}{x}  - \sum_{n=1}^{\infty} \frac{2x}{n^2\pi^2 - x^2}
\end{equation}
and hence
\begin{equation}\label{ecot}
\cot x 
\ge 
\frac{1}{x} - \sum_{n =1} ^{\infty}\frac{2x}{n^2 \pi^2 - \frac{\pi^2}{4}}
~=~ 
\frac{1}{x} - \frac{2x}{\pi^2} \sum_{ n=1}^{\infty} 
(\frac{1}{n-\frac{1}{2}} - \frac{1}{n+ \frac{1}{2}})
~=~ 
\frac{1}{x} - \frac{4x}{\pi^2}.
\end{equation}
Applying \eqref{ecot} with $x = \frac{m\pi}{k+1}$, which belongs to the
interval $]0, \frac{\pi}{6}]$, we obtain
\begin{equation}\label{an}
\cot (\frac{m\pi}{k+1})
\ge
 \frac{k+1}{m\pi} - \frac{4m}{(k+1)\pi},
\end{equation}
and we conclude by noting that 
$|q_m| ~=~ e^{-2\pi t_m} ~=~e^{-\pi  \cot \left( \frac{m \pi}{k+1} \right)}$.

\bigskip
\noindent{\bf Lemma 7.}$-$ 
{\em We have
$\frac{m+1}{m} \le e^{\frac{1}{m} - \frac{1}{2m(m+1)}}$.} 

\medskip
Indeed for any $x \in [0, 1[$, we have
\begin{equation}\label{elog}
- \log(1-x)
~=~ 
\sum_{n =1}^{\infty} \frac{x^n}{n} 
\le
x + \frac{1}{2}\sum_{ n =2}^{\infty} x^n
~=~ 
x + \frac{x^2}{2(1-x)}.
\end{equation}
Putting $x = \frac{1}{m+1}$ in \eqref{elog}, we obtain
\begin{equation}\label{elog1}
\log (\frac{m+1}{m})
\le 
\frac{1}{m+1} + \frac{1}{2m (m+1)}
~=~ 
\frac{1}{m} - \frac{1}{2m(m+1)}.
\end{equation}

\smallskip
From now on, we shall put
$u_n := \frac{n^k q_m^n}{(1- q_m^n)^2}$
for all $n \ge 1$. Hence we have $h_k(\frac{1}{2} + it_m)
= \sum_{n=1}^{\infty} u_n$ by formula \eqref{nf}.

\bigskip
\noindent{\bf Lemma 8.}$-$ 
{\em We have
$|u_{n+1}| < \frac{|u_n|}{2}$, for $n \ge m$. 
}

\medskip
Suppose that $n \ge m$. We have
\begin{eqnarray}\label{eqF}
\frac{|u_{n+1}|}{|u_n|} 
&=& 
(\frac{n+1}{n})^k ~|q_m| 
~\left(\frac{1-q_m^n}{1- q_m^{n+1}} \right)^2  \\
&\le&
(\frac{m+1}{m})^k ~|q_m|~ 
\left(\frac{1 + |q_m|}{1- |q_m|^2}\right)^2 \nonumber \\
&=&
(\frac{m+1}{m})^k ~\frac{|q_m|}{(1- |q_m|)^2}. \nonumber
\end{eqnarray}
Using the upper bounds of $|q_m|$ and of $\frac{m+1}{m}$ obtained
in Lemmas 6 and 7, and the inequality $k \ge 3m^2$, we get
\begin{equation}\label{eqF1}
(\frac{m+1}{m})^k |q_m|
\le
e^{\frac{4m}{k+1} - \frac{1}{m} - \frac{k}{2m(m+1)}}
\le
e^{\frac{4m}{3m^2+1} - \frac{1}{m} - \frac{3m}{2(m+1)}}
\le
e^{-3/4},
\end{equation}
where the last inequality follows from the fact that
$$
 \frac{4m}{3m^2 +1} - \frac{1}{m} - \frac{3m}{2(m+1)}  
 $$
 is equal to $-\frac{3}{4}$ 
when $m=1$, and is bounded by 
$$
\frac{4}{3m} -\frac{1}{m} -\frac{3m}{2(m+1)}   
~=~
 \frac{1}{3m} -\frac{3}{2} + \frac{3}{2(m+1)} 
 \le 
 \frac{1}{6} -\frac{3}{2} + \frac{1}{2}  
 = 
-\frac{5}{6} < - \frac{3}{4}
$$
when $m \ge 2$. 

Since $t_m \ge \frac{\sqrt{3}}{2}$, we have
$|q_m| = e^{-2\pi t_m} \le e^{-\pi \sqrt{3}}$, 
and we therefore deduce from \eqref{eqF}
and \eqref{eqF1} that 
\begin{equation}\label{eqF2}
\frac{|u_{n+1}|}{|u_n|} 
\le
\frac{e^{-3/4}}{(1- e^{-\pi \sqrt{3}})^2}
~<~ 
\frac{1}{2}.
\end{equation}

\medskip
\noindent{\bf Lemma 9.}$-$ 
{\em We have 
$|u_{n-1}| < \frac{|u_n|}{2}$ for $2 \le n \le m$. }

\medskip
Suppose that $2 \le n \le m$. We have
\begin{eqnarray}\label{eqF3}
\frac{|u_{n-1}|}{|u_n|} 
&=& 
(\frac{n-1}{n})^k ~|q_m|^{-1} 
\left(\frac{1 - q_m^n}{1- q_m^{n-1}} \right)^2 \\
&\le &
(\frac{m-1}{m})^k ~|q_m|^{-1} 
\left(\frac{1 + |q_m|^2}{1- |q_m|} \right)^2. \nonumber
\end{eqnarray}
Now we have $\log (\frac{m-1}{m}) \le - \frac{1}{m} - \frac{1}{2m^2}$
and $|q_m|^{-1} = e^{\pi \cot\left( \frac{m \pi}{k+1} \right)} \le e^{\frac{k+1}{m}}$,
from where it follows
\begin{equation}\label{eqF4}
(\frac{m-1}{m})^k ~|q_m|^{-1} 
\le
 e^{\frac{1}{m} - \frac{k}{2m^2}}
\le
e^{\frac{1}{m} - \frac{3}{2}} 
\le  e^{-1}
\end{equation}
as $k \ge 3m^2$ and $m \ge 2$.  We also have 
$|q_m| =  e^{-2\pi t_m} 
\le e^{-\pi \sqrt{3}} $. We deduce that
\begin{equation}\label{eqF5}
\frac{|u_{n-1}|}{|u_n|} 
\le
\frac{1}{e} \left(\frac{1 + |q_m|^2}{1 - |q_m|} \right)^2
\le 
\frac{1}{e}\left(\frac{1+ e^{-2\pi\sqrt{3}}}{1 - e^{-\pi \sqrt{3}}}\right)^2
~<~ 
\frac{1}{2}.
\end{equation}

\medskip
We now complete the proof of Proposition 1 when $k \ge 3m^2$.
The sign of $u_n$ is $(-1)^n$ for all $n \ge 1$.
We deduce from Lemma 8 that
$\frac{u_m}{2} + \sum_{n=m+1}^{\infty} u_n$
is an alternating series of real numbers whose sum is non-zero
and has sign $(-1)^m$. Similarly Lemma 9 
implies that the finite sum 
$\frac{u_m}{2} + \sum_{n =1}^{m-1} u_n$
is real, non-zero and has sign $(-1)^m$.
Hence $h_k(\frac{1}{2} + it_m) = \sum_{n= 1}^{\infty} u_n$
is a non-zero real number whose sign is $(-1)^m$.

\medskip
\subsection{A new expression of the function $h_k$ \nopunct}$~~$\\

\noindent{\bf Lemma 10.}$-$ 
{\em 
For any $z \in \mathfrak{H}$, we have
\begin{equation}\label{eH1}
h_k(z)=(2\pi)^{-k-1}k!\sum_{c=1}^{\infty}
\sum_{d\in \Z}c\Big({cz+d\over i}\Big)^{-k-1}.
\end{equation}
}

\medskip
For any $z \in \mathfrak{H}$, we have
\begin{equation}\label{eH2}
\sum_{d\in \Z} (z+d)^{-2}
~=~{\pi^2\over\sin^2(\pi z)}
~=~
(2\pi i)^2{q\over (1-q)^2}
~=~
(2\pi i)^2\sum_{n=1}^\infty n q^n,
\end{equation}
where $q = e^{2\pi i z}$.  By differentiating $k-1$ times
with respect to $z$, we obtain
\begin{equation}\label{eH3}
(-1)^{k-1}k!\sum_{d\in \Z} (z+d)^{-k-1}=(2\pi i)^{k+1}\sum_{n=1}^\infty n^kq^n.
\end{equation}
This equality, applied to $c z$ where $c$ is an integer $\ge 1$, 
can be written as
\begin{equation}\label{eH4}
(2\pi)^{-k-1}k! \sum_{d \in \Z}\Big({c z+d\over i}\Big)^{-k-1}
=
\sum_{n=1}^\infty n^k q^{n c}.
\end{equation}
Thus the right-hand side of \eqref{eH1} is equal to 
$\sum_{c=1}^{\infty}\sum_{n=1}^\infty cn^kq^{nc}$, i.e. to
$h_k(z)$ by formula \eqref{vf}.

\medskip
\subsection{Proof of Proposition 1 when $3m^2 \ge k +1$ \nopunct}$~~$\\

In this subsection, we assume that $m$ is an integer such that
$1 \le m \le \left[ \frac{k+1}{6} \right]$ and $3m^2 \ge k+1$.
This implies that we have $m \ge 2$ and $k \ge 11$.
As before, we let $t_m :=\frac{1}{2} \cot ( \frac{m\pi}{k+1} )$ 
and $z_m := \frac{1}{2} + i t_m$.

\smallskip
It follows from Lemma 10 that the sum
\begin{equation}\label{vexp}
\sum_{ c =1}^{\infty} 
\sum_{ d\in \Z} c\left(\frac{c z_m + d}{i}\right)^{-k-1}
\end{equation}
is real. We shall show that it is non-zero and that its sign
is $(-1)^m$. This, along with Lemma 10, will then
imply Proposition 1 in the case considered. 

\smallskip
Let $\theta_m :={m\pi\over k+1}$. Note that we have
\begin{equation}\label{eG1}
z_m ={1\over 2}(1+i\cot \theta_m)
= \frac{ie^{-i \theta_m}}{2 \sin \theta_m}
= i e^{-i\theta_m}|z_m|,
\end{equation}
 \begin{equation}\label{eG2}
z_m -1
={1\over 2}(-1+ i \cot \theta_m)
= \frac{i e^{ i \theta_m}}{2 \sin \theta_m}
= i e^{i\theta_m}|z_m|.
\end{equation}
The contribution to the sum \eqref{vexp} of the terms indexed by the
ordered pairs $(c,d)$ equal to $(1, 0)$ and $(1,-1)$ is therefore
\begin{equation}\label{eG3}
( \frac{z_m}{i})^{-k-1} 
+
(\frac{ z_m -1}{i})^{-k-1} 
=
2 \cos((k+1)\theta_m) |z_m|^{-k-1}
=
2(-1)^m |z_m|^{-k-1}.
\end{equation}
It is real and its sign is $(-1)^m$. To conclude, it will therefore
suffice to prove that the sum of the absolute 
values of the other terms of \eqref{vexp} is strictly less
than $2|z_m|^{-k-1}$.

\bigskip
\noindent{\bf Lemma 11.}$-$ 
{\em For any $x \in [0, \frac{1}{9}]$, we have
$\frac{1+ 9x}{1+x} \ge e^{\alpha x}$,
where $\alpha = 9\log\left(\frac{9}{5} \right)$.
}

\medskip
The functions
$x \mapsto x+1$ and $x \mapsto e^{\alpha x}$ are convex
and increasing in the interval $[0, \frac{1}{9}]$, hence
the function $x \mapsto (1+x)e^{\alpha x} - (1+9x)$ is convex
in this interval. Since it takes the value $0$ at both of its endpoints,
it is less than or equal to $0$ in the interval considered.

\bigskip
\noindent{\bf Lemma 12.}$-$ 
{\em  
For any integer $d \ge 1$ and 
any $x \in [0, \frac{1}{4d-1}]$, 
we have
\begin{equation}\label{eG7}
\frac{1+ (2d+1)^2x}{1+ (2d-1)^2x} 
\ge
\frac{1+ 9x}{1+x}
\phantom{m}\text{and}\phantom{m}
\frac{1+(2d+1)^2x}{1+x} 
\ge
\left( \frac{1+ 9x}{1+x}\right)^d.
\end{equation}  
}

The first inequality follows from the identity
\begin{equation}\label{eG4}
\frac{1+ (2d+1)^2x}{1+ (2d-1)^2x} 
~-~
\frac{1+ 9x}{1+x}
~=~
\frac{8(d-1)x(1- (4d-1)x)}{(1+x)(1+(2d-1)^2x)}.
\end{equation}

\smallskip
We prove the second inequality by induction on $d$.
It is clear when $d~=~1$. Suppose that $d \ge 2$.
If $x$ belongs to the interval $[0, \frac{1}{4d-1}]$,
it belongs to the interval $[0, \frac{1}{4d-5}]$
and hence we have 
${1 + (2d -1)^2x\over 1+x} \ge \big({1+9x\over 1+x}\big)^{d-1}$
by the induction hypothesis. By multiplying this
inequality by the first inequality of \eqref{eG7}, we obtain
the second inequality.

\medskip
For any integer $j \ge 0$, let $\A_j$ be the set of integers 
$d \ge 1$ such that $d+ \frac{1}{2}$ lies in 
the set $[jt_m, ~ (j+1)t_m [$. 
Its cardinality is bounded by 
\begin{equation}\label{eq-3.3.1}
t_m+ 1 
~=~ 
\frac{1}{2} \cot (\frac{m \pi }{k+1}) + 1
\le
\frac{k+1}{2 m \pi } + 1
\le
\frac{1}{2\pi} \sqrt{3(k+1)} + 1
\le
\sqrt{k+1}
\end{equation}
since $m \ge \sqrt{\frac{k+1}{3}}$ and $k +1 \ge 12$.

\bigskip
\noindent{\bf Lemma 13.}$-$ 
{\em  
We have
$\sum_{d \in \A_0} \frac{|z_m+ d|^{-k-1}}{|z_m|^{-k-1}} 
< \frac{1}{6000}$.  
}

\medskip
Let $d \in \A_0$. We have
\begin{equation}\label{eqG1}
\left( \frac{|z_m + d|}{|z_m|} \right)^2
~=~
\frac{(d+\frac{1}{2})^2 + t_m^2}{(\frac{1}{2})^2 + t_m^2}
~=~
\frac{1 + (2d+1)^2\tan^2 \theta_m}{1+ \tan^2\theta_m}.
\end{equation}
We have by hypothesis $d+ \frac{1}{2} < t_m$,
hence $\tan^2 \theta_m < \frac{1}{(2d+1)^2}$.
Thus $\tan^2\theta_m$ is bounded above by $\frac{1}{4d-1}$
and also by $\frac{1}{9}$. We now obtain from 
Lemmas 11 and 12 that 
\begin{equation}\label{eqG2}
\left( \frac{|z_m + d|}{|z_m|}\right)^2
\ge
\left(\frac{1 + 9 \tan^2\theta_m}{1+ \tan^2\theta_m} \right)^d
\ge
e^{d \alpha \tan^2 \theta_m}
\ge
e^{d \alpha \theta_m^2},
\end{equation}
where $\alpha = 9 \log(\frac{9}{5})$,
and consequently that
\begin{equation}\label{eqG3}
\left( \frac{|z_m + d|}{|z_m|}\right)^{-k-1}
\le
e^{-d\alpha \frac{k+1}{2}\theta_m^2}
~=~
e^{-\frac{d\alpha m^2 \pi^2 }{2(k+1)}}
\le
e^{-\frac{d \alpha \pi^2}{6}},
\end{equation}
since by hypothesis $3m^2 \ge k+1$.
This implies that
\begin{equation}\label{eqG5}
\sum_{d \in \A_0} \frac{|z_m+d|^{-k-1}}{|z_m|^{-k-1}} 
\le
\sum_{d=1}^{\infty} e^{-\frac{d \alpha \pi^2}{6}}
~=~
\frac{1}{e^{\frac{\alpha\pi^2}{6}} - 1}
~<~
\frac{1}{6000}.
\end{equation}

\bigskip
\noindent{\bf Lemma 14.}$-$ 
{\em We have
$\sum_{d =1}^{\infty} \frac{|z_m+ d|^{-k-1}}{|z_m|^{-k-1}} < \frac{5}{16}$ and 
$\sum_{d\!\!\!\le\!\!\!-2} \frac{|z_m+ d|^{-k-1}}{|z_m|^{-k-1}} <\frac{5}{16}$.
}

\medskip

If $d \in \A_j$, we have $|z_m+d|^2 
~=~ (d+ \frac{1}{2})^2 + t_m^2 \ge (j^2 + 1)t_m^2$.
Also since $t_m \ge \frac{\sqrt{3}}{2}$, we have
\begin{equation}\label{eq-3.3.4}
|z_m|^2 ~=~ \frac{1}{4} + t_m^2  \le \frac{4}{3}t_m^2. 
\end{equation}
It follows that $\frac{|z_m+d|^2}{|z_m|^2}$ is bounded below by 
$\frac{3(j^2+1)}{4}$, hence by $\frac{3}{2}$ when $j=1$ and
by $\frac{3j^2}{4}$ when $j \ge 2$. Consequently, using
\eqref{eq-3.3.1},
\begin{equation}\label{eq-3.3.2}
\sum_{d \in \A_1 } \frac{|z_m+ d|^{-k-1}}{|z_m|^{-k-1}}
\le
(\frac{3}{2})^{-(k+1)/2}~\sqrt{k+1}
\le
(\frac{3}{2})^{-6} ~\sqrt{12},
\end{equation}
where the last inequality follows from the inequality $k+1 \ge 12$
and the fact that the function $x \mapsto x(\frac{3}{2})^{-x}$ 
is decreasing for $x \ge 12$.  Similarly we have for $j \ge 2$,
\begin{equation}\label{eq-3.3.3}
\sum_{d \in \A_j } \frac{|z_m+ d|^{-k-1}}{|z_m|^{-k-1}}
\le
(\frac{3j^2}{4})^{-6} ~\sqrt{12}.
\end{equation}
Using Lemma 13 and the inequalities 
\eqref{eq-3.3.2} and \eqref{eq-3.3.3}, we obtain
\begin{equation}\label{eqG7}
\sum_{d = 1}^{\infty} \frac{|z_m+ d|^{-k-1}}{|z_m|^{-k-1}} 
\le
\frac{1}{6000} ~+~ (\frac{3}{2})^{-6} \sqrt{12}
~+~ (\frac{3}{4})^{-6} \sqrt{12} (\zeta(12)-1)
~<~ 
\frac{5}{16}.
\end{equation}
This proves the first assertion of Lemma 14. The second assertion is
deduced from the first one 
by noticing that $|z_m+d|=|z_m-d-1|$ for any integer $d \ge 1$.

\bigskip
\noindent{\bf Lemma 15.}$-$ 
{\em  Let $c$ be an integer $\ge 2$.  We have
\begin{equation}\label{eG9}
\sum_{d \in \Z} c \frac{|c z_m+ d|^{-k-1}}{|z_m|^{-k-1}} 
< c^{1-k}((\frac{\sqrt{3}}{2})^{-k-1} + 3 ).
\end{equation}  }

For any $r \in \Z$, let $\rC_r$ be the set of integers $d\in \Z$
for which the real part of $c z_m +d$ belongs to the
interval $[c(r - \frac{1}{2}), c(r + \frac{1}{2})[$. The cardinality
of $\rC_r$ is $c$. For all $d \in \rC_r$, using \eqref{eq-3.3.4},
we have
\begin{equation}\label{eG11}
c |c z_m+d|^{-k-1} 
\le
\begin{cases}
c^{-k} t_m^{-k-1}  
\le 
c^{-k} (|z_m| \frac{\sqrt{3}}{2})^{-k-1}   &\text{ if }  r=0,\\
 c^{-k} ||r| - \frac{1}{2} + it_m|^{-k-1} 
= 
c^{-k} |z_m +  |r| - 1|^{-k-1} &\text{ if } r \ne 0.
\end{cases}
\end{equation}
We then deduce from \eqref{eG11} and Lemma 14 that
\begin{eqnarray}\label{eG13}
\sum_{d \in \Z} c~ \frac{|c z_m +d|^{-k-1}}{|z_m|^{-k-1}}
\!\!\!\!\!\!\!\!\!\!
&&\le
c^{1-k} \left( (\frac{\sqrt{3}}{2})^{-k-1} ~+~
 2 ~\sum_{d = 0}^{\infty} \frac{|z_m + d|^{-k-1}}{|z_m|^{-k-1}}\right) \\
&&\le
 c^{1-k} \left( (\frac{\sqrt{3}}{2})^{-k-1} ~+~ 2(1+ \frac{5}{16}) \right). \nonumber
\end{eqnarray}
This completes the proof of Lemma 15.

\medskip
We deduce from Lemmas 14 and 15 that the sum
of the terms ${c|cz_m+d|^{-k-1}\over |z_m|^{-k-1}}$,
extended to all tuples $(c,d) \in \Z^2$ distinct from $(1, -1)$ 
and $(1,0)$ and for which $c \ge 1$,
is bounded above by 
\begin{eqnarray}\label{eqG11}
&&
\frac{5}{8} + \sum_{c\!\! \ge\!\! 2} c^{1-k} 
( (\frac{\sqrt{3}}{2})^{-k-1} + 3)
\le
\frac{5}{8} + \sum_{c\!\! \ge \!\! 2} c^{-10} 
( (\frac{\sqrt{3}}{2})^{-12} + 3)\\
&& 
\phantom{mmmmmmmmmmmmmm}
=~~
\frac{5}{8} + ( (\frac{\sqrt{3}}{2})^{-12} + 3) (\zeta(10) -1)
~<~ \frac{2}{3}. \nonumber
\end{eqnarray}
This completes the proof of Proposition 1, when $3m^2 \ge k+1$.

\medskip
\section{On the multiplicity of zeros of quasi-modular forms \\
with algebraic Fourier coefficients}

\medskip
\subsection{The algebra of quasi-modular forms for $\SL_2(\Z)$ \nopunct}$~~$ \\

The normalized Eisenstein series $\E_4$ and $\E_6$ have 
the Fourier expansions
\begin{eqnarray}\label{e4.1}
\E_4 (z) &=& 1 + 240\sum_{n=1}^{\infty} \sigma_3(n) q^n, \\
\E_6 (z) &=& 1 - 504\sum_{n=1}^{\infty} \sigma_5(n) q^n.
\end{eqnarray}
They are algebraically independent over $\C$ and the $\C$-algebra
of modular forms for $\SL_2(\Z)$ is none other than $\C[\E_4, \E_6]$.
We have $\E_4^3 - \E_6^2 = 1728 \Delta$, where the {\it modular
discriminant} $\Delta$ is defined by
\begin{equation}\label{e4.2}
\Delta(z) = q \prod_{n=1}^{\infty} (1 - q^n)^{24}.
\end{equation}

\smallskip
Formula \eqref{NEis} which defines Eisenstein series when 
$k \ge 4$ retains a meaning when $k=2$
and allows us to define a function $\E_2$ on the Poincar\'e half-plane, 
periodic with period $1$, by its Fourier expansion
\begin{equation}\label{e4.3}
\E_2(z) = 1 - 24\sum_{n=1}^{\infty} \sigma_1(n) q^n,
\end{equation}
where as usual $q = e^{2 \pi i z}$. We have $\E_2 = \frac{\D\Delta}{\Delta}$,
where $\D$ is the differential operator
$\frac{1}{2\pi i}\frac{d}{dz} = q \frac{d}{dq}$. It follows that $\E_2$
is not modular, but it satisfies for every element $({a\ b\atop c\ d})$
of $\SL_2(\Z)$ the relation
\begin{equation}\label{e4.4}
\E_2 ({az+b\over cz+d})=(cz+d)^{2}\E_2(z) + \frac{12c(cz+d)}{2\pi i}.
\end{equation}

\smallskip
The functions $\E_2, \E_4, \E_6$ are algebraically independent over $\C$.
We call {\it quasi-modular forms} for $\SL_2(\Z)$ the elements of the $\C$-algebra
$\C[\E_2, \E_4, \E_6]$ generated by these three functions. This algebra
is graded by the weight, agreeing that $\E_k$ is of weight $k$ for
$k \in \{2,4,6\}$.

\smallskip
Let $\X, \Y, \bZ$ be three indeterminates. For any polynomial ${\rm P} \in \C[\X, \Y, \bZ]$,
let $\psi_{\rm P} := {\rm P}(\E_2, \E_4, \E_6)$. The map $\rP \mapsto \psi_{\rP}$
is an isomorphism from the polynomial algebra $\C[\X, \Y, \bZ]$ to the algebra
$\C[\E_2, \E_4, \E_6]$ of quasi-modular forms for $\SL_2(\Z)$. The Fourier
coefficients of $\psi_{\rP}$ belong to a subfield ${\rm K}$ of $\C$ if and only
if $\rP$ belongs to ${\rm K}[\X, \Y, \bZ]$. The function $\psi_{\rP}$ is a
quasi-modular form of weight $k$ if and only if the polynomial 
$\rP$ is isobaric of weight $k$, when assigning $\X, \Y, \bZ$
the weights $2, 4, 6$ respectively.

\medskip
\subsection{Derivations of the algebra of quasi-modular forms \nopunct}$~~$\\

Ramanujan noticed that one has
\begin{equation}\label{rama}
\D \E_2 = \frac{1}{12}(\E_2^2 - \E_4), ~~~~
\D \E_4 = \frac{1}{3}(\E_2 \E_4  - \E_6), ~~~~
\D \E_6 = \frac{1}{2} (\E_2 \E_6 - \E_4^2).
\end{equation}
Consequently, $\D$ defines a graded derivation of degree $2$
of the graded algebra $\C[\E_2, \E_4, \E_6]$
of quasi-modular forms for $\SL_2(\Z)$.

\medskip
Let us denote by $\frac{\partial}{\partial \E_2}$ the derivation of the algebra
$\C[\E_2, \E_4, \E_6]$ which, under the isomorphism 
$\rP \mapsto \psi_{\rP}$, corresponds to the derivation $\frac{\partial}{\partial \X}$
of the algebra of polynomials $\C[\X, \Y, \bZ]$. It is graded of degree $-2$.

\medskip
The two derivations $\frac{\partial}{\partial \E_2}$ and $\D$ stabilise
$\K[\E_2, \E_4, \E_6]$ for any subfield $\K$ of~$\C$. Their bracket
$[\frac{\partial}{\partial \E_2}, \D] =
 \frac{\partial}{\partial \E_2} \D - \D\frac{\partial}{\partial \E_2}$  
 is the derivation of the algebra $\C[\E_2, \E_4, \E_6]$ which
 maps any quasi-modular form $\psi$ of weight $k$
to $\frac{k}{12}\psi$ : it suffices to check this when $\psi$ is
one of the functions $\E_2, \E_4, \E_6$, and it then follows
from formula \eqref{rama}.

\bigskip
\noindent{\bf Lemma 16.}$-$ 
{\em  Let $r$ be an integer $\ge 1$.  The bracket 
$[\frac{\partial}{\partial \E_2}, \D^r]$
maps any quasi-modular form $\psi$ of weight $k$ 
for $\SL_2(\Z)$ to $\frac{(k+r-1)r}{12} \D^{r-1}\psi$.  }

\medskip
We shall argue by induction on $r$. When $r=1$, our assertion
has already been proved. Suppose that $r\ge 2$. It
follows from the induction hypothesis that, for any quasi-modular form
$\psi$ of weight $k$ for $\SL_2(\Z)$, we have
\begin{eqnarray*}
[\frac{\partial}{\partial\E_2}, \D^r] \psi 
&=& 
[\frac{\partial}{\partial\E_2}, \D^{r-1}] \D\psi 
+ \D^{r-1} [\frac{\partial}{\partial \E_2}, \D] \psi \\
&=& 
\frac{(k+r)(r-1)}{12}\D^{r-1}\psi  + \frac{k}{12} \D^{r-1}\psi 
~=~ \frac{(k+r-1)r}{12} \D^{r-1}\psi
\end{eqnarray*}

\bigskip
\noindent{\bf Lemma 17.}$-$ 
{\em  Let $r$ and $j$ be integers such that $0\le j \le r$. For
any modular form $f$ of weight $k$ for $\SL_2(\Z)$, we have
\begin{equation}\label{brac1}
\left( \frac{\partial}{\partial \E_2}\right)^j \D^r\!f 
~=~
\left( \prod_{i=1}^{j} \frac{(k+r -i)(r-i+1)}{12} \right)\D^{r-j}\!f.
\end{equation}
Here, as usual, an empty product is considered to be equal to $1$.
 }

\medskip
We prove this by induction on $j$. When $j=0$, it is clear.
Suppose that $j\ge1$. By the induction hypothesis, we have
\begin{eqnarray*}
\left(\frac{\partial}{\partial\E_2} \right)^j\D^r\!f 
&=&
\frac{\partial}{\partial\E_2} 
\left( \left(\frac{\partial}{\partial\E_2} \right)^{j-1} \D^r\!f \right) \\
&=&
\left(\prod_{i=1}^{j-1} \frac{(k+r -i)(r-i+1)}{12} \right)
\frac{\partial}{\partial\E_2} \D^{r-j+1}\! f 
\end{eqnarray*}
As $f$ belongs to the algebra $\C[\E_4, \E_6]$ of modular forms, 
we have $\frac{\partial f}{\partial \E_2} =0$. Thus
$\frac{\partial}{\partial \E_2} \D^{r-j+1}\!f$ is equal to 
$[\frac{\partial}{\partial\E_2}, \D^{r-j +1}]f$, hence to
$\frac{(k+r-j)(r-j+1)}{12} \D^{r-j}\! f$ by
Lemma 16. This completes the proof of Lemma 17.

\medskip
\noindent
{\bf Remark 11.}$-$ Lemma 17 remains valid when we replace $f$ by the 
quasi-modular form $\E_2$ (with $k=2$ in this case). The proof is 
analogous to that of Lemma 17, except that we replace there the relation
$\frac{\partial f}{\partial \E_2} =0$ by $\frac{\partial \E_2}{\partial \E_2} =1$.

\medskip
\subsection{On the simplicity of zeros of certain quasi-modular forms \nopunct}$~~$\\

As in subsection 4.1, let $\psi_{\rP} = \rP(\E_2, \E_4, \E_6)$ for any polynomial
$\rP \in \C[\X, \Y, \bZ]$. We denote by $\overline{\Q}$ the algebraic closure
of $\Q$ in $\C$.

\medskip
\noindent{\bf Proposition 2.}$-$
{\em Let $\rP$ and $\rm Q$ be two irreducible
elements of the unique factorization domain $\overline{\Q} [\X,\Y, \bZ]$
which are not associates of each other. The functions
$\psi_{\rP}$ and $\psi_{\rm Q}$ have no common zeros in the
Poincar\'e upper half-plane $\mathfrak{H}$.
}

\medskip

Let us argue by contradiction, by assuming that there
exists one common zero. We denote it by $a$. The map 
$u: {\rm R} \mapsto \psi_{\rm R}(a)$ is a ring
homomorphism from $\overline{\Q}[\X, \Y, \bZ]$ to $\C$. 
Its kernel $\mathfrak{p}$ is a prime ideal of 
$\overline{\Q}[\X, \Y,  \bZ]$. It contains
$\rP$ and $\rm Q$ which are two irreducible elements 
of $\overline{\Q}[\X, \Y,  \bZ]$, not associates of each other.
The height of $\mathfrak{p}$ is therefore at least $2$. Consequently,
the ring $\overline{\Q}[\E_2(a),  \E_4(a), \E_6(a)]$
which is the image of $u$ has Krull dimension at most $1$. Since it
is a $\overline{\Q}$ algebra of finite type, its Krull
dimension is the transcendence degree of its fraction
field, i.e. of $\overline{\Q}(\E_2(a),  \E_4(a), \E_6(a))$.
Now this transcendence degree is at least $2$ by
a theorem of G. V. Chudnovsky \footnote[2]{\tiny{ The theorem 
of Chudnovsky was stated in a different manner
in terms of periods and quasi-periods of elliptic 
curves (see \cite{GC}, Chapter 7, Theorem 3.1); for the equivalence of
the two formulations, see for example in \cite{MW}, the
equivalence of Theorem 3 and its Corollary 1.}\par}.

\bigskip
\noindent
{\bf Corollary 1.}$-$
{\em Let $\rP$ and $\rm Q$ be two non zero elements 
of the unique factorization domain $\overline{\Q}[\X, \Y, \bZ]$
without common irreducible factors.The functions $\psi_{\rP}$ and
$\psi_{\rm Q}$ have no common zeros in $\mathfrak{H}$.
}

\bigskip
\noindent
{\bf Example 1.}$-$ 
Let $f$ be a quasi-modular form whose Fourier coefficients 
are algebraic over $\Q$. If $f$ vanishes at $i$, $f$ is a multiple of $\E_6$ 
in the ring of quasi-modular forms. If it vanishes at $e^{\pi i/3}$, it is a
multiple of $\E_4$.

\bigskip
\noindent
{\bf Corollary 2.}$-$
{\em Let $\rP$ be a non zero element of the unique factorization domain 
$\overline{\Q}[\X, \Y, \bZ]$. Let $a$ be a zero of $\psi_{\rP}$
in $\mathfrak{H}$ and $e$ be its order. There exists
an irreducible factor $\rm R$ of $\rP$, unique up to multiplication
by a scalar, such that $\psi_{\rm R}(a)=0$. The function $\psi_{\rm R}$ 
has a simple zero at $a$ and the $\rm R$-adic valuation of $\rP$
is $e$.
}

\smallskip
The existence of $\rm R$ is clear. Its uniqueness up to 
multiplication by a scalar follows from Proposition 2.
It follows from Proposition 2 that $e = v_{\rm R}(\rP) e'$,
where $e'$ is the order of $\psi_{\rm R}$ at $a$
and $v_{\rm R}(\rP)$ is the $\rm R$-adic valuation
of $\rP$. Recall that $\D$ denotes the differential operator
$\frac{1}{2\pi i}\frac{d}{dz} = q\frac{d}{dq}$.  
There exists a polynomial
${\rm Q} \in \overline{\Q}[\X, \Y, \bZ]$ such that
$\D\psi_{\rP} = \psi_{\rm Q}$. The order at
$a$ of $\D\psi_{\rP}$ is $e-1$. By Proposition 2, it
is equal to $v_{\rm R}({\rm Q})e'$. But then we have
$(v_{\rm R}(\rP) - v_{\rm R}({\rm Q}))e' =1$.  It implies that
$e'=1$. We then have $v_{\rm R}(\rP) = e$.

\bigskip
\noindent
{\bf Corollary 3.}$-$
{\em Let $\rP$  be an element of $\overline{\Q}[\X, \Y, \bZ]$
without multiple factors. The quasi-modular form
$\psi_{\rP}$ has only simple zeros in $\mathfrak{H}$.
}

\smallskip
This is an immediate consequence of Corollary 2.

\bigskip
\noindent
{\bf Corollary 4.}$-$
{\em The quasi-modular form $\E_2$ has only simple zeros in $\mathfrak{H}$.
}

\smallskip
This result follows from Corollary 3 by taking $\rP = \X$.

\bigskip
\medskip
Corollary 4 can also be deduced from the
first equality of  \eqref{rama}, as remarked by H. Saber and 
A. Sebbar in \cite{SS}, p. 1786. Much more precise 
information about the zeros of $\E_2$ can be found
in \cite{IJT} and \cite{WY}.

\bigskip
\subsection{Multiple zeros of iterated derivatives of modular forms \nopunct}$~~$\\

\medskip
\noindent
{\bf Theorem 6.}$-$
{\em Let $f$ be a non-zero modular form of weight $k>0$ for $\SL_2(\Z)$
whose Fourier coefficients are algebraic over $\Q$, and let $r$ be an integer 
$\ge 1$. If the $r$-th derivative of $f$ has a zero at a point $a$ in $\mathfrak{H}$ 
of order $e\ge 2$, then $f$ has a zero of order $e + r$ at this point.
}

\smallskip
By hypothesis, the function $\D^r\!f = (2\pi i)^{-r}\frac{\partial^r \!f}{\partial z^r}$
has a zero of order $e \ge 2$ at~$a$. Let $\rP$ be the element of 
$\overline{\Q}[\X, \Y, \bZ]$ such that $\psi_{\rP} = \D^r\!f$, 
with the notations of subsection 4.1.
By Corollary 2 of Proposition 2, there exists an irreducible factor
$\rm R$ of $\rP$ in $\overline{\Q}[\X, \Y, \bZ]$, unique up to multiplication
by a scalar, such that $\psi_{\rm R}(a) =0$. Further $\psi_{\rm R}$ has 
a simple zero at $a$ and we have $v_{\rm R}(\rP)=e$. 
We distinguish two cases :

\medskip
$a)$ {\it The polynomial $\rm R$ is of degree $\ge 1$ in $\X$.}

\smallskip
The polynomial $\rm R$ is then also irreducible when we consider it as
an element of the principal ideal domain $\overline{\Q}(\Y, \bZ)[\X]$,
and we have $v_{\rm R}(\frac{\partial {\rP}}{\partial \X}) = e - 1$.
Now we have
$$
\psi_{\frac{\partial{\rP}}{\partial \X}} 
= \frac{\partial}{\partial \E_2}\D^{r-1}\!f
= \frac{(k+r-1)r}{12} \D^r \!f
$$
by Lemma 17, hence $\D^{r-1}f$ has a zero of order $e-1$ at $a$.
It implies that $\D^r f$ has a zero of order $e-2$ at $a$, and this
contradicts our hypothesis.
Hence case~$a)$ cannot occur.

\medskip
$b)$ {\it The polynomial $\rm R$ belongs to $\overline{\Q}[\Y, \bZ]$.}

\smallskip
In this case, if we consider $\rP$ as a polynomial in $\X$ with
coefficients in $\overline{\Q}[\Y, \bZ]$, $\rm R$ divides the
coefficient of each power of $\X$ in $\rP$. Consequently, for
every integer $j\ge 0$,  $\rm R$ divides the polynomial 
$\frac{\partial^j \rP}{\partial \X^j}$, and therefore 
$\frac{\partial^j}{(\partial \E_2)^j}\D^r\!f$ vanishes at $a$.
It follows from Lemma 17 that $\D^{r-j}\!f(a) = 0$ for
$0 \le j \le r$. Thus $f$ and its iterated derivatives
of order $\le r$ vanish at $a$. Since the $r$-th derivative of
$f$ has a zero of order $e$ at $a$, $f$ has a zero of
order $e+r$ at $a$.

\bigskip
\subsection{Simplicity of the zeros of iterated derivatives 
of Eisenstein series \nopunct}$~~$

\bigskip
\noindent
{\bf Theorem 7.}$-$
{\em 
Let $k$ be an even integer $\ge 2$. For any integer $r \ge 1$, 
the $r$-th derivative of the function $\E_k$ has only simple zeros
in the Poincar\'e upper half-plane.
}

\smallskip
When $k\ge 4$, this follows from Theorem 6 applied to the modular form 
$f = \E_k$~: indeed, the zeros of $\E_k$ in $\mathfrak{H}$ have order 
$\le 2$ by subsection 1.2.

\smallskip
When $k=2$, Theorem 6 remains true for $f = \E_2$, although $\E_2$
is not a modular form : the proof is identical, by appealing to
Remark 11 of subsection~4.2. We conclude by noting that the
zeros of $\E_2$ in $\mathfrak{H}$ are simple, by Corollary 4 
of Proposition 2.

\medskip
\section{Zeros of $\E_k'$ in $\gamma\D$, for $\gamma \in \SL_2(\Z)$}$~~$

\subsection{Derivatives of modular functions \nopunct}$~~$\\

Let $k$ be an even integer $\ge 2$.
Let us denote by $\mM_k$ the vector space of meromorphic functions
$f$ on the Poincar\'e upper half-plane $\mathfrak{H}$ such that
\begin{equation}\label{e5.1}
f \left(\frac{az+b}{cz+d}\right) ~=~ (cz+d)^k f(z),
\end{equation}
holds in $\mathfrak{H}$, for any $\begin{pmatrix} a & b \\ c & d \end{pmatrix} \in \SL_2(\Z)$. 
This implies that
\begin{equation}\label{e5.2}
f'\left(\frac{az+b}{cz+d}\right) ~=~ (cz+d)^{k+2} f'(z) ~+~ kc (cz+d)^{k+1}f(z)
\end{equation}
holds in $\mathfrak{H}$, for any $\begin{pmatrix} a & b \\ c & d \end{pmatrix} \in \SL_2(\Z)$.

For every non-zero function $f \in \mM_k$, let $\varphi_f$ denote 
the meromorphic function on $\mathfrak{H}$ defined by
\begin{equation}\label{e5.3}
\varphi_f (z) := z + k \frac{f}{f'}(z).
\end{equation}
Since we have $f'(-\frac{1}{z}) = z^{k+2} f'(z) + k z^{k+1} f(z)$
by formula \eqref{e5.2}, we have 
\begin{equation}\label{e5.4}
\varphi_f (z) = z^{-k-1}\frac{f'(-\frac{1}{z})}{f'(z)}
\end{equation}
in $\mathfrak{H}$. To the best of our knowledge, the functions $\varphi_f$ have
been introduced as well as investigated in this context by 
H. Saber and A. Sebbar in \cite{SS}. These functions  have
the following important equivariance property : 

\bigskip
\noindent{\bf Lemma 18.}$-$ 
{\em Let $f$ be a non-zero element of ${\mM}_k$. We have
\begin{equation}\label{e5.5}
\varphi_f({az+b\over cz+d})
={a\varphi_f(z)+b\over c\varphi_f(z)+d}
\end{equation}
in $\mathfrak H$, for any $({a\ b\atop c\ d})\in \SL_2(\Z)$.}

Formula \eqref{e5.5} follows from formula \eqref{e5.3} when 
$({a\ b\atop c\ d})=({1\ 1\atop 0\ 1})$ and from formula \eqref{e5.4} when 
$({a\ b\atop c\ d})=({0\ -1\atop 1\ \ 0})$. These two matrices generate 
the group $\SL_2(\Z)$. The lemma follows. 

\bigskip
\noindent
{\bf Remark 12.}$-$ One can prove the identity \eqref{e5.5} 
by a direct computation :  it may be useful for example if one tries to 
generalize Lemma 18 when $\SL_2(\Z)$ is replaced by any 
Fuchsian group. Indeed, it follows from 
formulae \eqref{e5.1} and \eqref{e5.2} that the following equalities of 
meromorphic functions hold in $\mathfrak{H}$ :
 \begin{eqnarray*}
\varphi_f(\frac{az+b}{cz+d})
 &=&
 \frac{az+b}{cz+d} + \frac{k}{(cz+d)^2\frac{f'}{f}(z)+ kc(cz+d) } \\
 &=&
 {az+b\over cz+d}+{1\over (cz+d)({cz+d\over \varphi_f(z)-z}+c)} \\
& =&
{az+b\over cz+d}+{\varphi_f(z)-z\over (cz+d)(c\varphi_f(z)+d)}\\
&=&
{((az+b)c+1)\varphi_f(z) +(az+b)d-z\over (cz+d)(c\varphi_f(z)+d)}\\
& =&
{(acz+ad)\varphi_f(z) +bcz+bd\over (cz+d)(c\varphi_f(z)+d)} 
~=~
{a\varphi_f(z)+b\over c\varphi_f(z)+d}.
\end{eqnarray*}

\medskip
\subsection{\bf Zeros of derivatives of  modular  functions \nopunct}$~~$\\

In this subsection, $k$ denotes an even integer $\ge 2$ and $f$ a 
non-zero element of ${\mM}_k$. The zeros of the derivative $f'$ 
of $f$ in $\mathfrak{H}$ are of two different types : \par

$a)$ The  multiple zeros of $f$ : a zero of $f$ of multiplicity $e \ge 2$ 
is a zero of $f'$ of multiplicity $e-1$. We shall call 
them {\it the trivial  zeros of $f'$}. The set of these zeros is 
stable under the action of $\SL_2(\Z)$.\par

$b)$ The zeros of $f'$ which are not zeros of $f$. We shall 
call them {\it the non trivial  zeros of $f'$} and we shall 
denote the set of those zeros of $f'$ by $\bZ(f')$. If $\tau$ is an element of 
$\bZ(f')$ and $({a\ b\atop c\ d})$ an element of $\SL_2(\Z)$ with $c\not=0$, 
${a\tau+b\over c\tau+d}$ cannot be a zero of $f'$ by formula \eqref{e5.2}. \par

In particular, it follows from $b)$ that in the upper half plane 
$\mathfrak{H}$ a fixed point 
of a matrix $({a\ b\atop c\ d})~\in~\SL_2(\Z)$ with $c \ne 0$ cannot 
belong to $\bZ(f')$. In other words, $\bZ(f')$ does not meet the orbits
of $i$ and of $e^{\pi i/3}$ under the action of $\SL_2(\Z)$.

\bigskip
\noindent{\bf Lemma 19.}$-$ 
{\em The poles of $\varphi_f$ are the points of $\bZ(f')$. 
Moreover, the order of a point  $\tau\in \bZ(f')$ as a pole of $\varphi_f$ is 
also its order as a zero of $f'$.}

\smallskip
Let $\tau \in \mathfrak{H}$. We deduce from  formula \eqref{e5.3} 
the following : if $\tau$ is 
neither a pole of $f$ nor a zero of $f'$, $\varphi_f$ is holomorphic 
at $\tau$; if $\tau$ is  a pole of order $e$ of $f$, it is a pole of 
order $e+1$ of $f'$, hence $\varphi_f$ is holomorphic at $\tau$; if $\tau$ is a 
zero of order $e\ge 2$ of $f$, it is a  zero of order $e-1$ of $f'$, 
hence $\varphi_f$ is holomorphic at~$\tau$;  finally, if $\tau$ is not a  
zero of $f$ but is a zero of order $e \ge 1$ of $f'$, it is a pole 
of order $e$ of~$\varphi_f$. This proves the lemma.

\bigskip
\noindent{\bf Lemma 20.}$-$ 
{\em Let $\tau\in\mathfrak{H}$, and let
$\gamma=({a\ b\atop c\ d})$ be an element of $\SL_2(\Z)$ such that  
$c\not=0$. In order for $\gamma \tau={a\tau+b\over c\tau+d}$
to be a non trivial zero of $f'$, it is necessary and sufficient 
that $\tau$ is a  zero of $\varphi_f+{d\over c}$, and then the order 
of $\gamma \tau$ as a zero of $f'$ is the same as the order of $\tau$ as a
zero of $\varphi_f + {d\over c}$.}

\smallskip
Indeed,  for $\gamma \tau$ to be a non trivial zero of $f'$, it is necessary and 
sufficient that  it is a pole of $\varphi_f$ (Lemma 19), i.e. $\tau$ is a pole
of ${a\varphi_f+b\over c\varphi_f+d}$ (Lemma 18), i.e. a 
zero of $\varphi_f+{d\over c}$. Then the order of $\gamma\tau$ as a zero of $f'$
is equal to the order of $\gamma\tau$ as a pole of $\varphi_f$, i.e. of  $\tau$ as a 
pole of ${a\varphi_f+ b\over c\varphi_f + d}$, and also of $\tau$ as
a zero of $\varphi_f+ {d\over c}$.

\bigskip
\noindent
{\bf Remark 13.}$-$ One can also deduce Lemma 20 from the formula
\begin{equation}\label{e5.6}
(cz+d)^{-k-1}{f'({az+b\over cz+d})\over f'(z)}=cz+d+kc{f(z)\over f'(z)}=c\varphi_f(z)+d, 
\end{equation}
which follows from formula \eqref{e5.2}.

\smallskip
\subsection{Derivative of the function $\varphi_f$\nopunct}$~~$\\

In this subsection, $k$ denotes an even integer $\ge 2$ and 
$f$ a non-zero element of ${\mM}_k$. By differentiating $(63)$, we get
\begin{equation}\label{e5.7}
\varphi'_f(z)=k+1-k{f(z)f''(z)\over f'(z)^2}= {\F_f(z)\over f'(z)^2}
\end{equation}
in $\mathfrak{H}$, where 
\begin{equation}\label{e5.8}
\F_f(z)=(k+1)f'(z)^2-kf(z)f''(z).
\end{equation}

We have the following noteworthy fact :

\bigskip
\noindent{\bf Lemma 21.}$-$ 
{\em The function $\F_f$ is an element  of ${\mM}_{2k+4}$.}

\smallskip
The function $\F_f$ is meromorphic in 
$\mathfrak{H}$. For any $({a\ b\atop c\ d})\in \SL_2(\Z)$, we have
\begin{eqnarray*}
f({az+b\over cz+d}) &=& (cz+d)^kf(z),\\
f'({az+b\over cz+d})
&=&
(cz+d)^{k+2}f'(z)+kc(cz+d)^{k+1}f(z), \\
f''({az+b\over cz+d})
&=&
(cz+d)^{k+4}f''(z)+2(k+1)c(cz+d)^{k+3}f'(z) \\
&& \phantom{mmmm} +k(k+1)c^2(cz+d)^{k+2}f(z),\\
\end{eqnarray*}
in $\mathfrak{H}$, hence
\begin{equation}\label{e5.9}
\F_f({az+b\over cz+d})=(cz+d)^{2k+4}\F_f(z).
\end{equation}

\bigskip
\noindent
{\bf Remark 14.}$-$ If $f$ is a modular form of weight $k$ for $\SL_2(\Z)$, $\F_f$  
is a cuspform of weight $2k+4$ for $\SL_2(\Z)$ : this follows from
Lemma 21, as the constant terms of the Fourier expansions 
of $f'$ and $f''$ are both equal to $0$.

\smallskip
\subsection{Properties of $\varphi_f$ : the real case \nopunct}$~~$\\

In this subsection, $k$ denotes an even integer $\ge 2$. We shall say that 
an element  $f$ of ${\mM}_k$ is {\it real} if we have $f(it) \in \R \cup\{\infty\}$ 
for every $t>0$. It is equivalent to saying that  
\begin{equation}\label{e7.0}
f(-\overline{z})=\overline{f(z)}
\end{equation}
holds in $\mathfrak{H}$.

\bigskip
\noindent
{\bf Example 2.}$-$
If there exists a real number $\T \ge 0$ such that $f$ is holomorphic
in the half-plane $\{z~\in~\C~|~\Im(z)~>~\T \}$, then $f$ admits in this 
half-plane a Fourier expansion of the form 
$f(z)=\sum_{n\in \Z} a_n(f)q^n$, where $q=e^{2\pi i z}$,
and $f$ is a real element of ${\mM}_k$ if and only if
 all the coefficients $a_n(f)$ are real.
 
\bigskip
\noindent{\bf Proposition 3.}$-$
{\em Let $f$ be a non-zero real element of  ${\mM}_k$ and let $\tau$ be a 
point of $\mathfrak{H}$. \par

$a)$ If the real part of $\tau$ belongs to ${1\over 2}\Z$ and 
$\tau$ is not a pole of $\varphi_f$, then the real part of 
$\varphi_f(\tau)$ is equal to that of $\tau$.\par

$b)$ If the modulus of $\tau$ is  $1$, the modulus of $\varphi_f(\tau)$ is $1$;  in 
particular $\tau$ is not a pole of $\varphi_f$.}

\smallskip
Let $\sigma\in {1\over 2} \Z$. Since $f$ is periodic of period $1$, one deduces 
from \eqref{e7.0} that all the values taken by $f$ on $\sigma + i\R_+^*$ 
lie in $ \R\cup\{\infty\}$. Hence all the values taken by $f'$ on 
this half-line lie in $i \R\cup\{\infty\}$, and the same holds for the function 
$z \mapsto \varphi_f(z)-z$. This proves assertion $a)$. 

\smallskip
One deduces from \eqref{e7.0} that $f'(-\overline{z})=-\overline{f'({z})}$ 
holds in $\mathfrak{H}$, 
and therefore also $\varphi_f(-\overline{z})=-\overline{\varphi_f({z})}$. 
To prove $b)$, it suffices to consider the case where $\tau$ is not a 
pole of $\varphi_f$, as the general case then follows by continuity. 
Then one concludes by noticing that 
\begin{equation}\label{e7.1}
|\varphi_f(\tau)|^2=\varphi_f(\tau)\overline{\varphi_f(\tau)}
=-\varphi_f(\tau)\varphi_f(-\overline{\tau})
=-\varphi_f(\tau)\varphi_f(-{1\over\tau})=1
\end{equation}
by Lemma 18.

\bigskip
\noindent{\bf Corollary 1.}$-$ {\em Let $f$ be a non-zero real element of 
${\mM}_k$. Let  $\gamma=({a\ b\atop c\ d})$ be an element of $\SL_2(\Z)$ such that
$c\not=0$, $|c|\not=|d|$ and $|c|\not=2|d|$. \par

$a)$ The function $\varphi_f$ does not take the value $-{d\over c}$ at any point 
of the boundary of $\D$ in~$\mathfrak H$.\par

$b)$ The function $f'$ has no non trivial zeros  on the boundary 
of  $\gamma\D$ in $\mathfrak{H}$.}

\smallskip
According to our hypotheses, $-{d\over c}$ is distinct from $-1$, $1$, 
$-{1\over 2}$, ${1\over 2}$. Hence, by Proposition~3, $\varphi_f$ cannot 
take the value $-{d\over c}$ 
at a point of $\mathfrak{H}$ whose real part is ${1\over 2}$ or $-{1\over 2}$, 
or whose modulus is $1$. In particular, it cannot take the value  $-{d\over c}$ 
at a point of the boundary of $\D$ in~$\mathfrak{H}$. This proves assertion $a)$.\par

Assertion $b)$ follows from $a)$ by Lemma 20.

\medskip
\subsection{\bf Properties of $\varphi_f$ when $f$ is an Eisenstein series \nopunct}$~~$\\

In this subsection, $k$ denotes an even integer $\ge 4$ and $\E_k$ the 
normalized Eisenstein series of weight $k$ for $\SL_2(\Z)$. 
To simplify the notation, we shall denote by $\varphi_k$ the  meromorphic 
function $\varphi_{\E_k}$. We hence have
\begin{equation}\label{e7.2}
\varphi_k(z)
=
z + k{\E_k(z)\over\E'_k(z)}=z^{-k-1}{\E'_k(-{1\over z})\over\E'_k(z)}.
\end{equation}
The derivative of $\varphi_k$ is given, using 
formulae \eqref{e5.7} and \eqref{e5.8}, by 
\begin{equation}\label{e7.3}
\varphi'_k(z)={\F_k(z)\over\E'_k(z)^2},
\end{equation}
where
\begin{equation}\label{e7.4}
\F_k(z)=(k+1)\E'_k(z)^2-k\E_k(z)\E''_k(z).
\end{equation}
It follows from Remark 14 of subsection 5.3 that $\F_k$ is a 
cusp form of weight $2k+4$ for $\SL_2(\Z)$. It has a simple zero at the cusp $\infty$, 
as is seen from the Fourier expansion of the right hand-side of \eqref{e7.4}. 
The set of its zeros in $\mathfrak{H}$ is therefore stable under the 
action of $\SL_2(\Z)$ and the weighted number of its zeros 
modulo $\SL_2(\Z)$ (counted with their multiplicities, 
with weight ${1\over 2}$ for those in the orbit of $e^{\pi i/2}$,  
${1\over3}$ for those in the orbit of $e^{\pi i/3}$,  and $1$ for the others) 
is equal to ${2k+4\over12}-1={k-4\over 6}$.

\bigskip
\noindent $a)$ {\it Behaviour of $\varphi_k$ on the 
half-line ${1\over 2}+i[{\sqrt{3}\over 2},+\infty[$}

\smallskip
By Lemma 19, the poles of $\varphi_k$ lying on the 
half-line ${1\over 2}+i[{\sqrt{3}\over 2},+\infty[$ are the non trivial 
zeros of $\E'_k$ lying on this
half-line. It follows from  Theorems 1 and 2 that their number is 
$n=[{k-4\over 6}]$. Moreover they belong to the open half-line  
${1\over 2}+i]{\sqrt{3}\over 2},+\infty[$
and are simple zeros of $\E'_k$, hence simple 
poles of $\varphi_k$. Let us write them 
$$
{1\over 2}+ib_1, \dots, {1\over 2}+ib_n,
$$ 
with $b_1>\dots>b_n>{\sqrt{3}\over2}$, and let us 
fix the convention that $b_0:=+\infty$ and $b_{n+1}:={\sqrt{3}\over 2}$. 

For any real number  $t \ge {\sqrt{3}\over 2}$ distinct from $b_1,\cdots, b_n$, 
the complex number  
$\varphi_k({1\over2}+it)$ has  real part ${1\over 2}$ by Proposition 3, $a)$, 
hence can be written as
${1\over2}+iv_k(t)$, with $v_k(t)\in \R$. The  variations of the function $v_k$ 
are described in the following two  lemmas. 

\bigskip
\noindent{\bf Lemma 22.}$-$ 
{\em $a)$ We have $v_k({\sqrt{3}\over2})=-{\sqrt{3}\over2}$ when $k\equiv0\bmod 6$ and $v_k({\sqrt{3}\over2})={\sqrt{3}\over2}$ when $k\not\equiv0\bmod 6$.\par
$b)$ We have $\lim_{t\to+\infty}v_k(t)=(-1)^{k\over 2}\infty$.\par 
$c)$ We have $\lim_{t\to b_m^-}v_k(t)=(-1)^{{k\over 2}+m}\infty$ for $1\le m \le n$.\par
$d)$ We have $\lim_{t\to b_m^+}v_k(t)=(-1)^{{k\over 2}+m-1}\infty$ for $1\le m \le n$.}

\smallskip
It follows from Proposition 3, $b)$ that $\varphi_k({1\over 2}+i{\sqrt{3}\over 2})$ 
has modulus $1$. Hence $v_k({\sqrt{3}\over2})$ must be equal 
to either ${\sqrt{3}\over2}$ or  $-{\sqrt{3}\over2}$. It follows from formula 
$\eqref{e7.2}$ that we have  $v_k({\sqrt{3}\over2})={\sqrt{3}\over2}$ if 
and only if ${\E_k\over \E'_k}$ vanishes at the point 
${1\over 2}+i{\sqrt{3}\over 2}$; according to subsection 1.2, 
this is the case if and only if  $k\not\equiv0\bmod 6$. This proves assertion  $a)$.

One deduces from formula $(3)$ that when $t$ tends to $+\infty$, 
$\E_k({1\over 2}+it)$ tends to $1$ and 
$\E'_k({1\over 2}+it)$ is equivalent to ${4\pi ik\over \B_k}e^{-2\pi t}$. 
It then follows from $(72)$ that $v_k(t)$ is equivalent to 
$-{\B_k\over 4\pi}e^{2\pi t}$, hence tends to 
$(-1)^{k\over 2}\infty$ when $t$ tends to $+\infty$. 
This proves assertion~$b)$.

For $t>{\sqrt{3}\over 2}$, we have $\E_k({1\over 2}+it)\in \R_+^*$ (subsection 1.2, 
Remark~2) and $i\E'_k({1\over 2}+it)\in \R$. Moreover the function 
$t\mapsto i \E'_k({1\over 2}+it)$ vanishes only at the points $b_1, \cdots, b_n$, 
and has simple zeros at these points. It therefore  keeps a constant sign on 
each interval $]b_{m+1},b_{m}[$, where $0 \le m \le n$ (with $b_0=+\infty$ 
and $b_{n+1}={\sqrt{3}\over2}$), and these 
signs are alternatively positive and negative. 
This sign is $(-1)^{k\over2}$ on $]b_{1},b_{0}[$ by the previous discussion, 
and is therefore  $(-1)^{{k\over 2}+m}$ on $]b_{m+1},b_{m}[$ for $0\le m \le n$. 
It follows that for $1\le m \le n$, the  imaginary part of  ${\E_k\over \E'_k}({1\over 2}+it)$ 
tends to $(-1)^{{k\over 2}+m}\infty$ when $t<b_m$ tends to 
$b_m$, and to $(-1)^{{k\over 2}+m-1}\infty$ when $t>b_m$ tends to $b_m$.
The same then holds for $v_k(t)$ by formula $\eqref{e7.2}$. 
This proves assertions $c)$ and $d)$.

\bigskip
\noindent{\bf Lemma 23.}$-$ 
{\em $a)$ In each interval $]b_{m},b_{m-1}[$, where $1\le m \le n$,
there is a unique point $c_m$ at which the derivative  $v'_k$ of $v_k$
vanishes. Moreover, $v'_k$ does not vanish at any point of 
the interval $]{\sqrt{3}\over2}, b_n[$. \par
$b)$ The zeros of  $\varphi'_k$ in $\mathfrak{H}$ are all simple. They are the 
elements of the orbits under the action of $\SL_2(\Z)$ of the points 
${1\over 2}+ i c_m$,  where  $1\le m \le n$,
together with those of the orbit of $e^{\pi i/3}$ 
when $k \equiv 0 \bmod 6$.}\par

Let $m$ be an integer such that $1\le m \le n$.  The function $v_k$ is differentiable 
(and even real analytic) in the interval $]b_{m},b_{m-1}[$, and it has the same infinite 
limit at both endpoints of this interval by Lemma 22. Its derivative 
$v'_k$ must therefore vanish at at least one point $c_m$ of this interval. 
This proves the existence of $c_m$ in assertion  $a)$.

The derivative $\varphi'_k$ of $\varphi_k$ then vanishes 
at the point ${1\over 2}+i c_m$. Hence the modular function $\F_k$ also
vanishes at this point by formula $\eqref{e7.3}$. We have thus produced  
$n=[{k-4\over 6}]$ zeros of $\F_k$ which 
are pairwise inequivalent modulo $\SL_2(\Z)$, and also not equivalent to 
$e^{\pi i/2}$ or to $e^{\pi i/3}$. Since the weighted number of zeros of $\F_k$ 
modulo $\SL_2(\Z)$ (counted with their multiplicities, with weight  ${1\over 2}$ 
for those in the orbit of $e^{\pi i/2}$,  ${1\over3}$ for those in the orbit 
of  $e^{\pi i/3}$,  and $1$ for the others) is equal to ${k-4\over 6}$, the zeros of 
$\F_k$ found previously are all simple, and there are no other zeros of $\F_k$ 
modulo $\SL_2(\Z)$, except the points of the orbit of  $e^{\pi i/3}$ which are  
simple zeros of $\F_k$ when $k\equiv 0\bmod 6$ and zeros of multiplicity~$2$ 
of $\F_k$ when $k\equiv 2\bmod 6$. When $k\equiv 0\bmod 6$, these latter points 
are simple zeros of $\varphi'_k$; when $k\equiv 2\bmod 6$, they are simple 
zeros of $\E'_k$ by subsection 1.2 and hence they are not 
zeros of $\varphi'_k$ by formula $\eqref{e7.3}$. 
We have thus proved assertion $b)$. 

The uniqueness in $a)$ as well as the last assertion of $a)$ follow from $b)$.

\bigskip
\noindent
{\bf Remark 15.}$-$
Assertion $a)$ of Lemme 23 implies that the function $v_k$ is 
strictly monotonic in the intervals $]b_m,c_m]$ and 
$[c_m,b_{m-1}[$, for $1\le m\le n$, as well as in the interval
$[{\sqrt{3}\over2}, b_n[$. We shall prove in Proposition 4 of 
subsection 5.6 that $\varphi_k$ 
does not take the value ${1\over 2}$ in $\D$. It follows that
$v_k$ does not vanish at any point of the intervals $]b_m, b_{m-1}[$, 
where $1\le m \le n+1$, and therefore
keeps a  constant sign on each of them. This sign is 
$(-1)^{\frac{k}{2} + m-1}$ by Lemma 22.

\bigskip
\noindent $b)$ {\it Behaviour of $\varphi_k$ on the half-line 
$-{1\over 2}+i[{\sqrt{3}\over 2},+\infty[$}

\smallskip
We have $\varphi_k(z+1)=\varphi_k(z)+1$ in $\mathfrak{H}$. Hence, 
with the notations of $a)$, the poles of $\varphi_k$ lying on the half-line 
$-{1\over 2}+i[{\sqrt{3}\over 2},+\infty[$ are 
$-{1\over 2}+ib_1, \dots, -{1\over 2}+ib_n$, and they are simple. 
Moreover, for any real number   $t\ge {\sqrt{3}\over 2}$ distinct 
from $b_1,\ldots,b_n$, we have
\begin{equation}\label{e7.5}
\varphi_k(-{1\over2}+it)
=
\varphi_k({1\over2}+ it)-1 = -{1\over2}+ i v_k(t),
\end{equation}
and the variations of $v_k$ have been studied in $a)$.

\bigskip
\noindent $c)$ {\it Behaviour of the function 
$\theta\mapsto\varphi_k(e^{i\theta})$ in the interval $[{\pi\over 3},{2\pi\over 3}]$}

\smallskip
All values of the function $\theta\mapsto \varphi_k(e^{i\theta})$ in 
the interval $[{\pi\over 3},{2\pi\over 3}]$ are of modulus~$1$ by 
Proposition 3, $b)$ of subsection 5.4. Its value at ${\pi\over 3}$ is 
equal to $e^{-{\pi i /3}}$ when $k\equiv0\bmod 6$ and to  
$e^{\pi i /3}$ when $k\not\equiv0\bmod 6$ by Lemma 22, $a)$.

Hence there exists a unique continuous  function  
$w_k:[{\pi\over 3},{2\pi\over 3}]\to \R$ such that
\begin{equation}\label{e7.6}
\varphi_k(e^{i\theta})=e^{iw_k(\theta)}
\end{equation}
for $\theta\in[{\pi\over 3},{2\pi\over 3}]$ and
\begin{equation}\label{e7.7}
w_k({\pi\over 3})
=
\begin{cases}
- {\pi\over3} &\text{ when }  k \equiv 0 \bmod 6,\\
{\pi\over3} &\text{ when }  k \not\equiv 0\bmod 6.
\end{cases}
\end{equation}

\bigskip
\noindent{\bf Lemma 24.}$-$ 
{\em The function $w_k$ is strictly increasing in 
$[{\pi\over 3},{2\pi\over 3}]$ and we have
\begin{equation}\label{e7.8}
w_k({2\pi\over 3})
=
\begin{cases}
{(k-2)\pi\over 3} & \text{ when } k\equiv 0\bmod 6,  \\
{k\pi\over 3}  &  \text{ when } k\equiv 2\bmod 6, \\
{(k+4)\pi\over 3} & \text{ when } k\equiv 4\bmod 6.
\end{cases}
\end{equation}
} \par

The function $w_k$ is real analytic in the interval 
$[{\pi\over 3}, {2\pi\over 3}]$. We get 
by differentiating $\eqref{e7.6}$
\begin{equation}\label{e7.9}
e^{i\theta}\varphi'_k(e^{i\theta})=w'_k(\theta)e^{iw_k(\theta)}
\end{equation}
for $\theta\in [{\pi\over 3},{2\pi\over 3}]$. It follows from 
Lemma 23, $b)$  that $\varphi'_k(e^{i\theta})$ does not vanish at 
any point of the open 
interval $]{\pi\over 3},{2\pi\over 3}[$. Therefore the 
function $w'_k$ does not vanish at any point of this interval, and keeps 
a constant sign on this interval. It follows that the function $w_k$ is strictly monotonic 
in the closed interval $[{\pi\over 3},{2\pi\over 3}]$.

\smallskip
Let us denote by $\rC$ the open circular arc consisting of the 
points  $e^{i\theta}$, where $\theta\in]{\pi\over3},{2\pi\over 3}[$. 
The function $\E'_k$ does not vanish at any  point of $\rC$
by Lemma 1 of subsection 2.2. 
For every  $\tau\in \rC$, we have
\begin{equation}\label{e8.0}
\varphi_k(\tau)=\tau^{-k-1}{\E'_k(-{1\over \tau})\over \E'_k(\tau)}
= \tau^{-k-1}{\E'_k(-\overline{\tau})\over \E'_k(\tau)}
=-\tau^{-k-1}{\overline{\E'_k(\tau)}\over \E'_k(\tau)}
\end{equation}
by  formulae $\eqref{e7.2}$ and $\eqref{eq1}$.

\smallskip
By definition (see subsection 2.1), the variation of the 
argument  of the function $\theta \mapsto\varphi_k(e^{i\theta})$ along the interval 
$[{\pi\over 3},{2\pi\over 3}]$ is $w_k({2\pi\over 3})-w_k({\pi\over 3})$. It follows 
from formula $\eqref{e8.0}$ that it is also equal to
$-{(k+1)\pi\over 3}-2 \A$ where, as in subsection 2.4, $\A$ denotes the 
variation of the argument of $\theta\mapsto \E'_k(e^{i\theta})$ 
along  $[{\pi\over 3},{2\pi\over 3}]$
when $k\not\equiv2\bmod 6$, and the limit when 
$\eta>0$ tends to $0$ of the variation of the argument of 
$\theta\mapsto \E'_k(e^{i\theta})$ along 
$[{\pi\over 3}+\eta,{2\pi\over 3}-\eta]$ when  $k\equiv2\bmod 6$.

\smallskip
The value of $\A$ has been computed in subsections 2.6, 2.7 
and 2.8. It is equal 
to $-{(k+2)\pi\over 3}$ when $k\equiv 4 \bmod 6$ and to $-{k\pi\over 3}$ 
when $k\not\equiv 4\bmod 6$. We have therefore
\begin{equation}\label{e8.1}
w_k({2\pi\over 3})-w_k({\pi\over 3})
=
\begin{cases}
-{(k+1)\pi\over 3}+2{(k+2)\pi\over 3}
={(k+3)\pi\over 3} & \text{ when } k\equiv 4\bmod 6, \\
-{(k+1)\pi\over 3} + 2{k\pi\over 3}
= {(k-1)\pi\over 3} &\text{ when } k\not\equiv 4\bmod 6.
\end{cases}
\end{equation}
One notes that $w_k({2\pi\over 3})>w_k({\pi\over 3})$ in both cases. 
The strictly monotonic  function $w_k$ is therefore strictly increasing 
in $[{\pi\over 3},{2\pi\over 3}]$. 
Finally, formula $\eqref{e7.8}$ follows from formulae 
$\eqref{e7.7}$ and $\eqref{e8.1}$.

\bigskip
\noindent
{\bf Remark 16.}$-$
As $\varphi_k(-\frac{1}{z})\varphi_k(z)=-1$, the function 
$\theta \mapsto w_k(\pi -\theta) + w_k(\theta)$ is constant on the
interval $[\frac{\pi}{3}, \frac{2\pi}{3}]$. By formula \eqref{e7.7}
and \eqref{e7.8}, its value is $\frac{(k-3)\pi}{3}$ when 
$k \equiv 0 \bmod{6}$, $\frac{(k+1)\pi}{3}$ when $k \equiv 2\bmod{6}$
and $\frac{(k+5)\pi}{3}$ when $k \equiv 4 \bmod{6}$.

\medskip
\subsection{\bf Zeros of $\E'_k$ in $\gamma\D$, for $\gamma \in \SL_2(\Z)$
\nopunct}$~~$\\

In this subsection, $k$ denotes an even integer $\ge 4$, $\E_k$  the 
normalized Eisenstein series of weight   $k$ for $\SL_2(\Z)$ and $\D$  
the closure of the  standard fundamental domain of  $\mathfrak{H}$ 
modulo $\SL_2(\Z)$. We denote by $\gamma$ a given element 
$\left({a\ b\atop c\ d}\right)$ of $\SL_2(\Z)$. We are interested in the 
number of zeros of $\E'_k$ lying in $\gamma\D$. 

\smallskip
We already know that all the zeros of $\E'_k$ in $\mathfrak{H}$ are
 simple (Theorem 7 of subsection 4.5). Counting those lying in 
 $\gamma\D$ with or without multiplicities therefore leads to the same result. 
 As in subsections 1.4 and 5.2, we shall call {\it trivial zeros} of $\E'_k$ those zeros 
of  $\E'_k$ which are also zeros of $\E_k$. Such zeros 
exist only when $k\equiv 2\bmod 6$ and they are then the elements of the orbit of 
$e^{\pi i/3}$ under the action of $\SL_2(\Z)$;  the set $\gamma\D$ hence 
contains  $2$ trival zeros of $\E'_k$ when $k\equiv 2\bmod 6$, namely 
$\gamma e^{\pi i/3}$ and $\gamma e^{2\pi i/3}$, and no trivial zeros otherwise.

\smallskip
As in subsections 1.4 and 5.2, we shall denote by $\bZ(\E'_k)$ the set of non 
trivial zeros of $\E'_k$ in $\mathfrak{H}$. We still have to count the number of points 
of the set  $\bZ(\E'_k)\cap\gamma\D$. By subsection 5.2,  $\bZ(\E'_k)$ 
does not meet the orbits of $i$ and of $e^{\pi i/3}$ under the 
action of $\SL_2(\Z)$. There is therefore
no need to assign weights to the points of $\bZ(\E'_k)$ while
counting them. 

\smallskip
When $\gamma=\left({1\ 0\atop 0\ 1}\right)$,  $\bZ(\E'_k)\cap\gamma\D$ 
consists of $[{k-4\over6}]$ points with real part ${1\over 2}$ and the same 
number with real part  $-{1\over2}$ (subsection 1.3, Theorems 1 and 2). 
More generally,  when $c=0$,  the cardinality of $\bZ(\E'_k)\cap\gamma\D$ is 
$2[{k-4\over6}]$, since  $\bZ(\E'_k)$ is translation invariant by $\Z$. 

\smallskip
We shall therefore suppose from now on that $c\not=0$. 
Let $\tau\in\D$; in order for $\gamma \tau$ to belong to $\bZ(\E'_k)$,  it is necessary 
and sufficient that we have $\varphi_k(\tau)=-{d\over c}$  (Lemma 20). The  
cardinality of $\bZ(\E'_k)\cap\gamma\D$ is therefore the number of 
roots of the equation $\varphi_k(z)=-{d\over c}$ in $\D$. Moreover, all these 
roots are simple ({\it loc. cit.}). So Theorem 5 of subsection 1.4 is a consequence 
of the following proposition :

\bigskip
\noindent{\bf Proposition 4.}$-$
{\em Let $\lambda$ be a rational number. \par
$a)$ When $|\lambda|<1$, the set $\varphi_k^{-1}(\lambda)\cap\D$ is empty.\par
$b)$ When $|\lambda|\ge1$, the set $\varphi_k^{-1}(\lambda)\cap\D$ has 
cardinality $[{k+2\over6}]$. It is contained in the interior of $\D$
when $|\lambda|\not=1$, and in $\rC=\{e^{i\theta}~|~\theta\in]{\pi\over 3},{2\pi\over3}[\}$ 
when $|\lambda|=1$.\par}

\smallskip
We can always write $\lambda$ in the from $-{d\over c}$ for some matrix
$\left({a\ b\atop c\ d}\right)\in\SL_2(\Z)$, with $c\not=0$. By the previous discussion,  
the roots of the equation $\varphi_k(z)=\lambda$ in $\D$ are simple. The poles
of $\varphi_k$ in $\D$ are, with the notations of subsection 5.5, the points 
of the form $\pm{1\over 2}+i b_m$, where $1\le m\le n$ 
and $n=[{k-4\over6}]$.

\bigskip
\noindent 1) {\it The case where $|\lambda|$ is distinct from $1$ and ${1\over 2}$}

\smallskip
In this case, the set $\varphi_k^{-1}(\lambda)\cap\D$ in contained in the 
interior of $\D$ (subsection 5.4, Corollary 1 of Proposition 3). 
When $z \in \D$ and $\Im(z)$ tends to $+\infty$, $\E_k(z)$ tends to $1$ 
and $\E'_k(z)$ is equivalent to $-{4\pi i k \over \B_k}e^{2\pi i z}$ by formula 
$(3)$, hence $\varphi_k(z)$ is equivalent to ${i\B_k\over 4\pi}e^{-2\pi i z}$ by 
formula $\eqref{e7.2}$. So there exists a real  number 
$\T>1$, such that $\T > b_1$ if $n \ge 1$ (i.e. $k \ge 10$), and such 
that all the points of $\varphi_k^{-1}(\lambda)\cap\D$ have
imaginary part $<\T$. For $\varepsilon >0$, let $\D_{\T,\varepsilon}$
denote the set of points in $\D$ with imaginary part $\le \T$ 
and also at a  distance $\ge \varepsilon$ from the poles  
$\pm{1\over 2} + i b_m$, where $1\le m\le n$. It follows 
from the residue
theorem that, when $\varepsilon$ is sufficiently small, we have
\begin{equation}\label{e8.2}
{\rm Card}(\varphi_k^{-1}(\lambda)\cap\D)
={1\over 2\pi i}\int_{\partial \D_{\T,\varepsilon}}{\varphi'_k(z)\over\varphi_k(z)-\lambda}\,dz
\end{equation}
where $\partial \D_{\T, \varepsilon}$ is the oriented boundary 
of $\D_{\T, \varepsilon}$ (see Figure 5). 

\medskip

\centerline{\includegraphics[width=9.6cm, height=8.4cm]{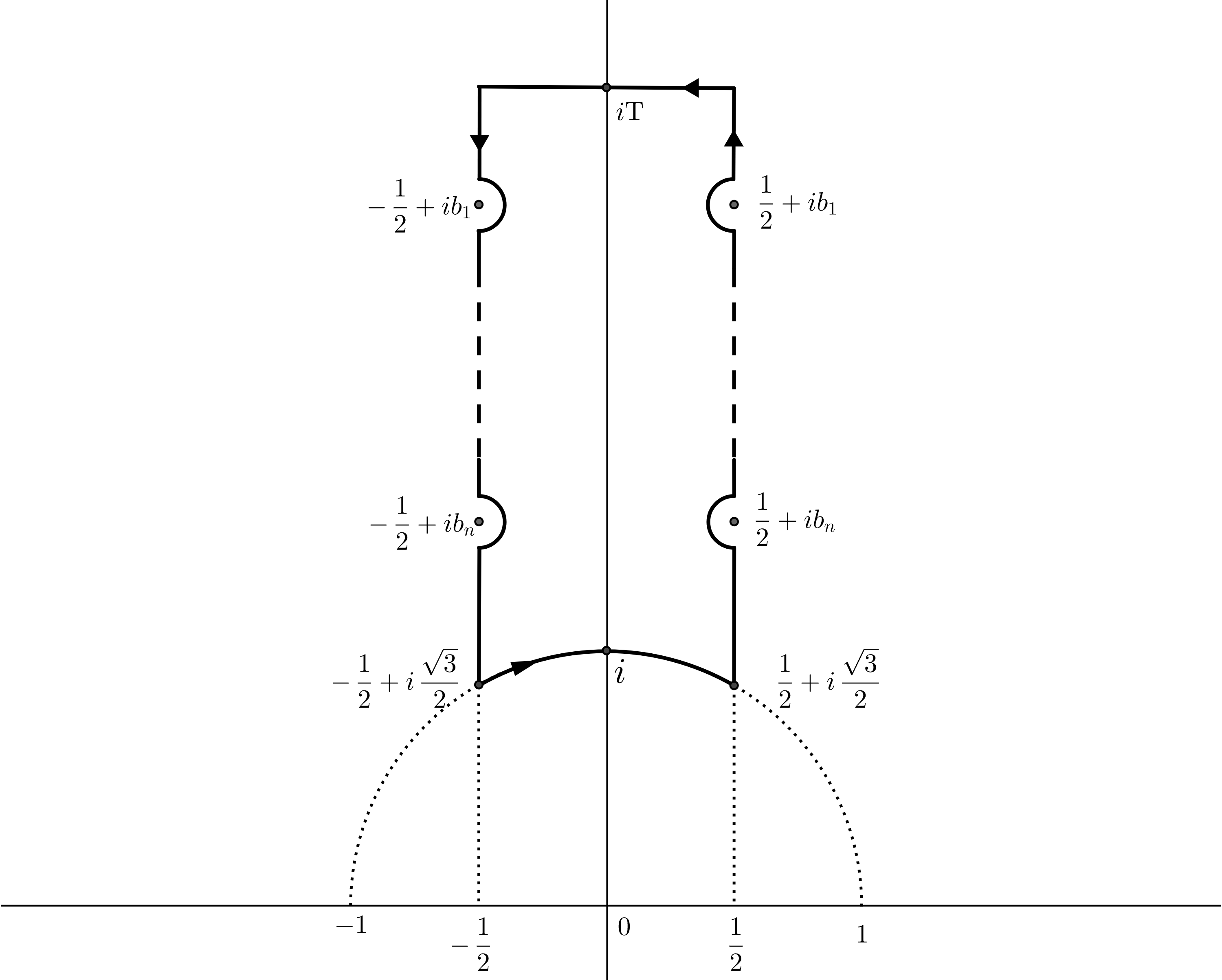}}

\centerline{{\it Figure 5.}  The contour of integration $\partial \D_{\T,\varepsilon}$.}

\bigskip
Formula $\eqref{e8.2}$ may also be written as 
\begin{equation}\label{e8.3}
{\rm Card}(\varphi_k^{-1}(\lambda)\cap\D)
={1\over 2\pi} (\text{variation of the argument of } \varphi_k-\lambda
\text{ along } \partial\D_{\T, \varepsilon} ).
\end{equation}

To evaluate the right-hand side, which does not depend 
on $\T$ and $\varepsilon$ (when $\T$ is sufficiently large 
and $\varepsilon$ sufficiently small), 
we shall decompose the contour of integration $\partial \D_{\T, \varepsilon}$ 
into finitely many parts, study the variation of the argument of 
$\varphi_k-\lambda$ along each of them, then take the limit of this 
variation when successively $\T$ tends to $+\infty$ 
and $\varepsilon$ to~$0$, and finally sum up these various  limits.

\medskip
\noindent $a)$ {\it Contribution of the upper horizontal  segment}

\smallskip
When $z\in \D$ and $\Im(z)$ tends to $+\infty$, we have seen that 
$\varphi_k(z)$ is  equivalent to ${i\B_k\over 4\pi}e^{-2\pi i z}$, and the same
 then holds for $\varphi_k(z)-\lambda$. It follows that the variation 
 of the argument of $\varphi_k - \lambda$ along  the oriented horizontal 
 segment joining ${1\over 2}+i \T$ to $-{1\over 2}+ i \T$ tends 
 to $2\pi$ when $\T$ tends to $+\infty$. 
 
\medskip
\noindent $b)$ {\it Contributions of the upper vertical segments when $n\ge 1$}

\smallskip
Suppose $n\ge 1$, i.e. $k\ge 10$. When $t$ varies in the interval 
$[b_1+ \varepsilon, \T]$, $\varphi_k({1\over 2}+it)$ is moving on the line 
${1\over 2}+ i \R$, starting  from the point ${1\over 2} + i v_k(b_1+\varepsilon)$
and ending at the point ${1\over 2}+iv_k(\T)$, where  $v_k$ is the function considered 
in subsection 5.5. By Lemma 22, $v_k(\T)$ when $\T$ tends to 
$+\infty$, and  $v_k(b_1+\varepsilon)$ when $\varepsilon$ tends to~$0$, have the
same infinite limit $(-1)^{\frac{k}{2}} \infty$. It follows that the  variation of the 
argument of $\varphi_k - \lambda$ 
along the oriented segment joining
${1\over2}+i(b_1+\varepsilon)$ to ${1\over2}+i \T$ tends to $0$ when $\T$ tends to 
$+\infty$ and then $\varepsilon$ tends to $0$.

The same holds for the variation of the argument of $\varphi_k-\lambda$ along the 
oriented vertical segment joining
$-{1\over2}+i \T$ to $-{1\over2}+i(b_1+\varepsilon)$.

\medskip
\noindent $c)$ {\it Contributions of the vertical segments lying between two consecutive poles}

\smallskip
Let $m$ be an integer such that $1\le m\le n-1$. A reasoning similar to that 
of $b)$ shows that the variation of the argument of $\varphi_k-\lambda$ 
along the oriented  vertical segment joining ${1\over2}+i(b_{m+1}+\varepsilon)$ 
to ${1\over2}+i(b_{m}-\varepsilon)$ tends to $0$ when $\varepsilon$
tends to~$0$. 

\smallskip
The same holds for the variation of the argument of 
$\varphi_k-\lambda$ along the oriented vertical segment joining 
$-{1\over2}+i(b_{m}-\varepsilon)$ to $-{1\over2}+i(b_{m+1}+\varepsilon)$.

\medskip
\noindent $d)$ {\it Contributions of the semicircles centered at the poles}

\smallskip
Let $m$ be an integer such that $1\le m\le n$. The function $\varphi_k$ has a  
simple pole at the point $z_m={1\over 2}+ib_m$, hence 
there exists $\alpha\in \C^*$ such that $\varphi_k(z)$ is  equivalent to 
$\alpha(z-z_m)^{-1}$ when $z$ tends to $z_m$, and the same then 
holds for $\varphi_k(z)-\lambda$. Therefore, the variation of the  
argument of $\varphi_k-\lambda$ along the semicircle of 
radius $\varepsilon$ centered at $z_m$, joining $z_m-i\varepsilon$ to 
$z_m+i\varepsilon$ and clockwise oriented,
tends to $\pi$ when $\varepsilon$ tends to $0$.

\smallskip
Similarly,  the variation of the   argument of $\varphi_k-\lambda$ along the 
semicircle of radius $\varepsilon$ centered at $z'_m=-{1\over 2}+ib_m$,
 joining $z'_m+i\varepsilon$ to $z'_m-i\varepsilon$  and clockwise oriented, 
 tends to $\pi$ when $\varepsilon$ tends to $0$.

\medskip
\noindent $e)$ {\it Remaining contribution when $n\ge 1$}
\smallskip
 
Suppose $n\ge 1$, i.e. $k\ge 10$. The remaining contribution concerns 
the part $\delta_{\varepsilon}$ of $\partial \D_{\T,\varepsilon}$ obtained
by moving successively along the oriented vertical segment joining 
$-{1\over 2}+i(b_n-\varepsilon)$ to $-{1\over 2}+i{\sqrt{3}\over 2}$, 
then along the clockwise oriented  arc of the unit circle
joining $e^{2\pi i/3}$ to $e^{\pi i/3}$ and finally along the oriented 
vertical segment joining ${1\over 2}+i{\sqrt{3}\over 2}$ to
${1\over 2}+i(b_n-\varepsilon)$.

\smallskip
When $z$ moves along $\delta_{\varepsilon}$, $\varphi_k(z)$ is 
successively moving :\par
--- on the line $-{1\over 2}+i \R$, starting from the point
$-{1\over 2}+iv_k(b_n-\varepsilon)$ and ending at the point 
$-{1\over 2}+ i v_k({\sqrt{3}\over 2})$;\par
--- on the unit circle, starting from the point $-{1\over 2}+iv_k({\sqrt{3}\over 2})$
and ending at the point ${1\over 2}+iv_k({\sqrt{3}\over 2})$, the  variation of the 
argument of $\varphi_k$ during this move being equal to 
$w_k({\pi\over 3})-w_k({2\pi\over 3})$, where $w_k$ is the 
function considered in Lemma~24;\par
--- on the line ${1\over 2}+i \R$, starting from the 
point ${1\over 2}+iv_k({\sqrt{3}\over 2})$
and ending at the point   ${1\over 2}+iv_k(b_n-\varepsilon)$.

\medskip
According to Lemma 22, c), the limit when $\varepsilon$ tends to 
$0$ of $v_k(b_n- \varepsilon)$ is
$(-1)^{{k\over2} + n}\infty=(-1)^{{k\over2}+[{k - 4\over 6}]}\infty$. By Lemma 22 
and formula $\eqref{e8.1}$, we have the following table :

\bigskip
\begin{center}
\begin{tabular}{|c|c|c|c|c|}
\hline
& $\lim\limits_{\varepsilon\to 0} v_k(b_n-\varepsilon)$ & $v_k({\sqrt{3}\over 2})$
&$w_k({\pi\over 3})-w_k({2\pi\over 3})$\\
\hline
&&&\\
$k\equiv0\bmod 6$ & $-\infty$ & $-{\sqrt{3}\over 2}$ &  $-{(k-1)\pi\over3}$ \\
&&&\\
\hline
&&&\\
$k\equiv2\bmod 6$&$+\infty$&${\sqrt{3}\over 2}$&$-{(k-1)\pi\over3}$ \\
&&&\\
\hline
&&&\\
$k\equiv4\bmod 6$&$+\infty$&${\sqrt{3}\over 2}$&$-{(k+3)\pi\over3}$ \\
&&&\\
\hline 
\end{tabular}
\end{center}

\bigskip
As an example, the trajectory of $\varphi_k(z)$ when $z$ 
varies along $\delta_{\varepsilon}$
is represented in figure 6 below, in the particular case 
where $k =12$ (and hence $n=1$). Double arrows in 
this picture indicate arcs on which $\varphi_k(z)$ passes
twice.

\bigskip

\centerline{\includegraphics[width=7.6cm, height=8.4cm]{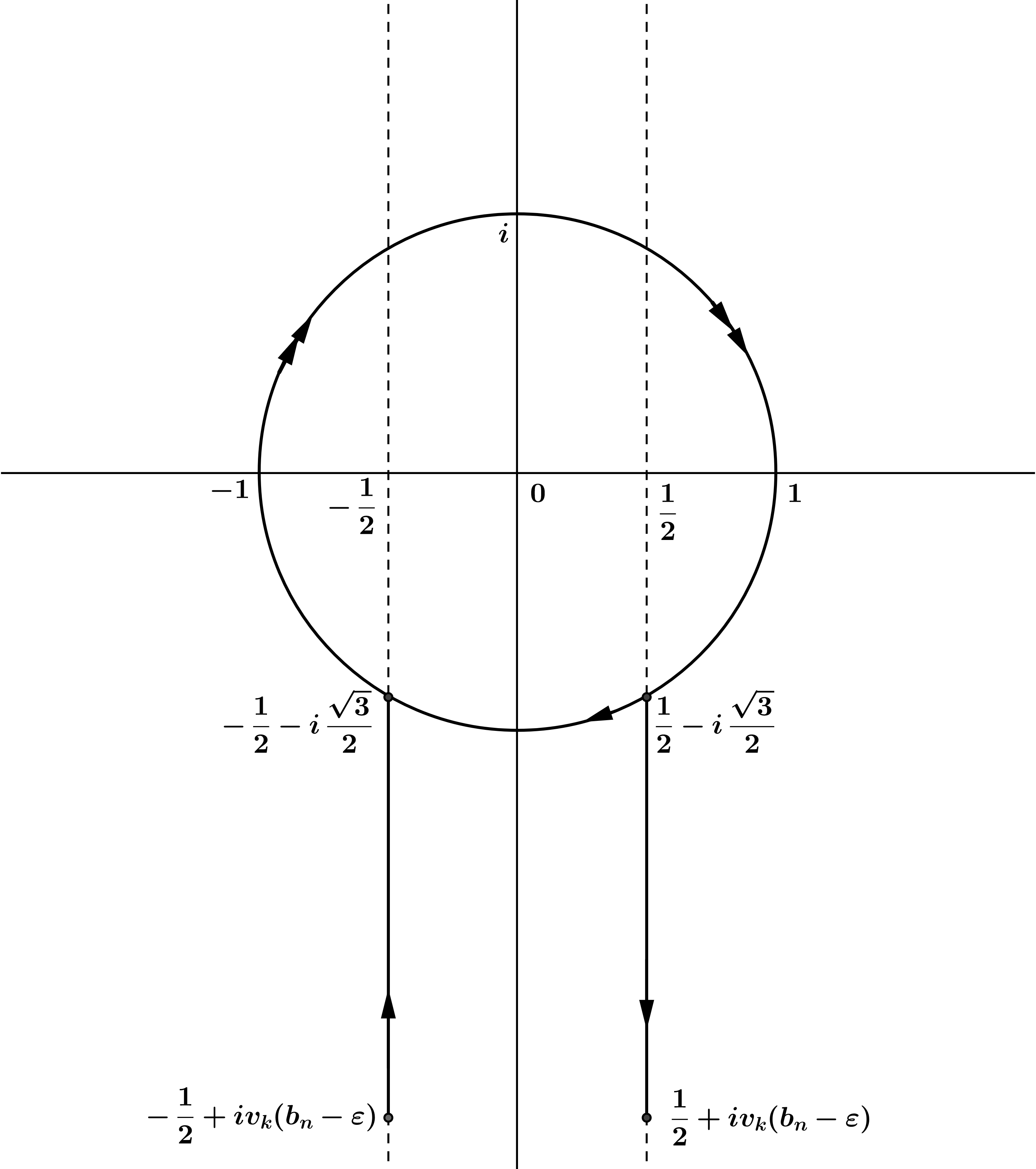}}

\medskip
\centerline{{\it Figure 6.}  Trajectory of $\varphi_k(z)$ when $z$ varies
along $\delta_{\varepsilon}$, for $k=12$.}

\bigskip
One deduces from the previous description of the trajectory of 
$\varphi_k(z)$ when $z$ varies along $\delta_{\varepsilon}$
and from the above table that :\par
$(i)$ when $|\lambda|>1$, the variation of the argument of 
$\varphi_k-\lambda$ along  $\delta_{\varepsilon}$ tends to~$0$
when $\varepsilon$ tends to $0$;\par
$(ii)$ when $|\lambda|<1$ and $|\lambda|\not={1\over 2}$, the limit 
when $\varepsilon$ tends to $0$ of the variation of the argument 
of $\varphi_k-\lambda$ along  $\delta_{\varepsilon}$ equals
$-{\pi\over3}-{(k-1)\pi\over3}=-{k\pi\over3}$ when $k\equiv0\bmod6$,
${\pi\over3}-{(k-1)\pi\over3}=-{(k-2)\pi\over3}$ when $k\equiv2\bmod6$ 
and ${\pi\over3}-{(k+3)\pi\over3}=-{(k+2)\pi\over3}$ 
when $k\equiv4\bmod6$, i.e. to
$-2[{k+2\over6}]\pi$ in all these three cases.

\medskip
\noindent $e')$ {\it Remaining contribution when $n=0$}
\smallskip

Suppose $n=0$, i.e. $k\le 8$. The remaining contribution 
then concerns the part $\delta_{\T}$ of $\partial \D_{\T,\varepsilon}$
analog to that considered in $e)$, except that $b_n -\varepsilon$ is 
now replaced by $\T$. One then proves in the same 
way as in $e)$ that the limit, 
when $\T$ tends to $+\infty$, 
of the variation of the argument of 
$\varphi_k-\lambda$ along  $\delta_{\T}$ is $0$ 
when $|\lambda|>1$, and is $-2[{k+2\over6}]\pi$
when $|\lambda|<1$ and $|\lambda|\not={1\over 2}\,$.

\medskip
\noindent $f)$ {\it Summing up the various contributions}
\smallskip

By summing up the limits of the various previously described 
contributions, one gets the total variation of the 
argument of $\varphi_k-\lambda$ 
along  $\partial \D_{\T,\varepsilon}$. It is equal to 
$2\pi+2n\pi=2[{k+2\over6}]\pi$ when $|\lambda|>1$, and to  
$2\pi+2n\pi-2[{k+2\over6}]\pi=0$ when $|\lambda|<1$ 
and $|\lambda|\not={1\over 2}$.
It then follows from $\eqref{e8.3}$ that we have
\begin{equation}\label{e8.4}
{\rm Card}(\varphi_k^{-1}(\lambda)\cap\D)
=
\begin{cases}
[{k+2\over6}] & \text{ when } |\lambda|>1, \\
0 & \text{ when } |\lambda|<1 \text{ and } |\lambda|\not={1\over 2}\,.
\end{cases}
\end{equation}

\bigskip
\noindent
{\bf Remark 17.}$-$
To get \eqref{e8.4}, key ingredients were Lemma 22 
and formula \eqref{e8.1}. On the other hand, we did neither
need the additional properties of $v_k$ stated
in Lemma 23 and in Remark 15, nor the fact that $w_k$ is an 
increasing function as stated  in Lemma 24.

\bigskip
\noindent 2) {\it The case where $|\lambda|={1\over 2}$ }
\smallskip

If the equation $\varphi_k(z)= {1\over 2}$ had a solution
in the interior of $\D$, then for each rational number $\lambda$ 
sufficiently close to $\frac{1}{2}$, the equation $\varphi_k(z)=\lambda$ would 
also have a solution, by the open mapping
theorem.This would contradict \eqref{e8.4}. 

\smallskip
If the equation $\varphi_k(z)={1\over 2}$ had a solution 
lying on the boundary of $\D$, this solution would lie on the open half-line
${1\over 2}+i]{\sqrt{3}\over2},+\infty[$ by Proposition 3 
of subsection 5.4.  In that case, for each rational number $\lambda$ sufficiently 
close to ${1\over 2}$, the equation $\varphi_k(z)=\lambda$ would 
have a solution in $\D\cup(1+ \D)$, by the open mapping
theorem, and hence either 
the equation $\varphi_k(z)=\lambda$ or the equation
$\varphi_k(z)=\lambda-1$ would have a solution in $\D$. 
This again would contradict \eqref{e8.4}.

\smallskip
We have thus proved  that the set $\varphi_k^{-1}({1\over2})\cap\D$ is empty. 
One proves in the same way that the set 
$\varphi_k^{-1}(-{1\over2})\cap\D$ 
is empty, or one deduces it from the identity 
$\varphi_k(-\overline{z})=-\overline{\varphi_k(z)}$.

\bigskip
\noindent 3) {\it The case where $|\lambda|=1$} 
\smallskip

If the equation $\varphi_k(z)=1$ had a solution in the interior of  $\D$, 
then for each rational number  $\lambda<1$ sufficiently close to 
$1$, the equation $\varphi_k(z)=\lambda$ would also have a solution,
 by the open mapping theorem. This would contradict \eqref{e8.4}. 
The set $\varphi_k^{-1}(1)\cap\D$ is therefore contained in 
the boundary of $\D$. It follows from Proposition 3, $a)$ of 
subsection 5.4 that it is in fact
contained in the open circular arc $\rC$ consisting of the points $e^{i\theta}$, 
where $\theta\in]{\pi\over3},{2\pi\over 3}[$.

By formula $(76)$, we have $\varphi_k(e^{i\theta})=e^{iw_k(\theta)}$ 
for every $\theta\in]{\pi\over3},{2\pi\over 3}[$. The cardinality of 
$\varphi_k^{-1}(1)\cap\D$ is therefore the number of elements  
$\theta\in]{\pi\over3},{2\pi\over 3}[$ for which we have
$w_k(\theta)\in 2\pi \Z$. 
It follows from Lemma 24 and  formula $(77)$ that this number is
equal to ${k\over 6}$ when $k\equiv 0\bmod 6$, to ${k-2\over 6}$ when 
$k\equiv 2\bmod 6$ and to ${k+2\over 6}$ when 
$k\equiv 4\bmod 6$, hence to $[{k+2\over 6}]$ in all these three cases.

One proves similarly, or one deduces from the  identity 
$\varphi_k(-\overline{z})= - \overline{\varphi_k(z)}$, that  
$\varphi_k^{-1}(-1)\cap\D$ is contained in $\rC$
and has $[{k+2\over 6}]$ elements.

\medskip
\subsection{The real locus of $\varphi_k$ \nopunct}$~~$\\

In this subsection, we keep the notations of subsections 5.5 and 5.6 : $k$ is an 
even integer $\ge 4$,  $\E_k$ is the normalized Eisenstein series of weight $k$ 
for $\SL_2(\Z)$ and $\varphi_k$ is the meromorphic function 
defined by $\varphi_k(z) = z + k {\E_k(z)\over \E'_k(z) }$ on the Poincar\'e
upper half-plane $\mathfrak{H}$. We shall present  
a pictorial interpretation of most of the results stated in subsections 5.5 and 5.6.

For this purpose, let us consider {\it the real locus} of the function $\varphi_k$, 
i.e. the set
$\varphi_k^{-1}( {\bP}_1(\R))$. It is stable under the action of $\SL_2(\Z)$, since
we have
\begin{equation}\label{eq85}
\varphi_k({az+b\over cz+d})={a\varphi_k(z)+b\over c\varphi_k(z)+d}
\end{equation}
in $\mathfrak{H}$ for every $({a\ b\atop c\ d})\in \SL_2(\Z)$ (Lemma 18). 
The study of this real locus can therefore be restricted to the  set
\begin{equation}\label{eq86}
{\rm R}_k=\varphi_k^{-1}({\bP}_1(\R))\cap\D,
\end{equation}
where $\D$ is the closure of the standard fundamental 
domain of $\mathfrak{H}$ modulo $\SL_2(\Z)$. 
As an example, ${\rm R}_k$ is represented for $k=16$ in 
Figure 7 below. 

\bigskip

\centerline{\includegraphics[width=9.6cm, height=8.4cm]{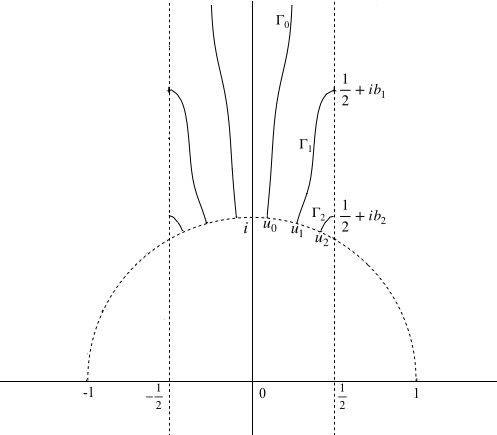}}

\medskip
\centerline{{\it Figure 7.}  The set $\rR_k=\varphi_k^{-1}({\bP}_1(\R))\cap\D$ 
for $k=16$.}

\bigskip

As $\varphi_k(-\overline{z})=-\overline{\varphi_k(z)}$ in $\mathfrak{H}$, the set 
$\rR_k$ is symmetric with respect to the imaginary axis.
So we can restrict its study to $\rR_k^+= \rR_k\cap \D^+$, 
where $\D^+$  is the set of points in $\D$ with non-negative real part. 
The set $\rR_k^+$ is described by the following proposition, 
which yields a  ``visual" explanation of Proposition~4~: 

\bigskip
\noindent{\bf Proposition 5.}$-$
{\em $a)$ The set $\rR_k^+$ does not meet the imaginary axis.\par
$b)$ The set $\rR_k^+$ has $n=[{k-4\over 6}]$ points with 
real part ${1\over 2}$. They are the poles of $\varphi_k$ located 
on the half-line ${1\over 2}+i[{\sqrt{3}\over 2},+\infty[$, 
i.e. the zeros of $\E'_k$ located on the open half-line 
 ${1\over 2}+i]{\sqrt{3}\over 2},+\infty[$. 
 They are therefore the points denoted by ${1\over 2}+i b_1,\ldots,{1\over 2}+ib_n$ 
 in subsection 5.5, $a)$, with $b_1>\ldots>b_n>{\sqrt{3}\over 2}$.\par
 
$c)$ The set $\rR_k^+$ has  $n+1=[{k+2\over 6}]$ points of modulus $1$. 
They are the points of $\D_+$ where $\varphi_k$ takes the value $1$ or $-1$. 
Let us denote them by $u_0, \ldots, u_n$, where 
$u_j=e^{i\alpha_j}$ and ${\pi\over2}>\alpha_0>\ldots>\alpha_n>{\pi\over 3}$.\par

$d)$ We have $\varphi_k(u_j)=(-1)^{{k\over 2}+j+1}$.\par

$e)$ The connected component $\Gamma_0$ of $u_0$ in $\rR_k^+$ is a 
real analytic curve, of which $u_0$ is an endpoint
and which has an infinite branch  asymptotic to the vertical line with 
abscissa~${1\over 4}$. The map $\varphi_k$ induces an analytic isomorphism from
$\Gamma_0$ onto the interval $[1,+\infty[$ when ${k\over 2}$ is odd, 
and onto the interval $]-\infty, -1]$ when~${k\over 2}$ is even.\par

$f)$ Let $j$ be an integer such that $1\le j\le n$. 
The connected component $\Gamma_j$ of $u_j$ in $\rR_k^+$  is a real analytic curve,
of which one endpoint is $u_j$ and the other one is  ${1\over 2}+i b_j$. The map 
$\varphi_k$ induces an analytic isomorphism from
$\Gamma_j$ onto the interval  $[1,+\infty]$ when ${k\over 2}+j$ is odd, 
and onto the interval  $[-\infty, -1]$ when ${k\over 2}+j$ is even, these intervals
being identified with their canonical images in $\bP_1(\R)$.\par 

$g)$ The curve $\Gamma_j$ is orthogonal to the unit circle at the point $u_j$, 
for $0 \le j \le n$. It is orthogonal to the vertical line with abscissa ${1\over 2}$ 
at the point ${1\over 2}+i b_j$, for $1\le j \le n$.\par

$h)$ The curves $\Gamma_j$, for $0\le j \le n$, are pairwise disjoint and
their union is~$\rR_k^+$.\par}

\smallskip
By differentiating formula $(3)$, we see that the function 
$\E'_k$ takes purely imaginary values and has no zeros on 
the half-line $i]0,+\infty[$. The same then holds
for the function $\varphi_k$ by formula $\eqref{e7.2}$. 
This proves assertion $a)$.

\smallskip
Let  $\tau$ be a point on the half-line ${1\over 2}+i[{\sqrt{3}\over2},+\infty[$, 
distinct from the poles ${1\over 2}+ib_1,\ldots,{1\over 2}+ib_n$ of $\varphi_k$. Then
$\varphi_k(\tau)$ has real part ${1\over 2}$ by Proposition 3,~$a)$, and is 
distinct from ${1\over 2}$ by Proposition 4, $a)$, hence is not real.
This proves assertion~$b)$. 

\smallskip
Let $\tau$ be a point of modulus $1$ of $\rR_k$. By Proposition 3, $b)$, 
$\varphi_k(\tau)$ has modulus~$1$; as it is a real number, it is 
equal to $1$ or $-1$. Conversely, by Proposition~4,~$b)$, any point of $\D$
 at which $\varphi_k$ takes the 
value $1$ or $-1$ has modulus $1$, and there are  $2(n+1)$ such points. 
Their set is symmetric with respect to the imaginary axis since 
$\varphi(-\overline{z})=-\overline{\varphi(z)}$, hence $n+1$ among them 
belong to~$\D^+$. This proves assertion $c)$.

\smallskip
By  formula $\eqref{e7.6}$, we have $\varphi_k(e^{i\theta})=e^{iw_k(\theta)}$ for 
$\theta\in[{\pi\over 3},{2\pi\over 3}]$, where $w_k$ is the function defined 
in subsection 5.5, $c)$. It follows that  $\alpha_0,\ldots,\alpha_n$ are 
the points of the interval $[{\pi\over 3},{\pi\over 2}]$, in descending  
order, at which the function $w_k$ takes a value in
$\pi\Z$. Since the function $w_k$ is strictly increasing (Lemma 24), the real 
numbers $\varphi_k(u_j)=e^{iw_k(\alpha_j)}$ are alternatively equal to $1$ and $-1$. 
It therefore suffices to prove assertion $d)$ when $j=0$. By the above,  
$w_k(\alpha_0)$ is the largest element of $\pi\Z$ which is smaller 
that $w_k({\pi\over 2})$. It follows from Remark 16 that we have
\begin{equation}\label{eq87}
w_k({\pi\over2})
=
\begin{cases}
{(k-3)\pi\over 6} & \text{ when } k\equiv 0 \bmod 6, \\
{(k+1)\pi\over 6} & \text{ when } k\equiv 2 \bmod 6, \\
{(k+5)\pi\over 6} & \text{ when } k\equiv 4 \bmod 6,
\end{cases}
\end{equation}
and therefore
\begin{equation}\label{eq88}
w_k(\alpha_0)
=
\begin{cases}
{(k-6)\pi\over 6} & \text{ when } k  \equiv 0 \bmod 6,\\
{(k-2)\pi\over 6}  & \text{ when } k\equiv 2 \bmod 6, \\
{(k+2)\pi\over 6}& \text{ when } k \equiv4 \bmod 6 \,.
\end{cases}
\end{equation}
In each of these three cases, we have 
$\varphi_k(u_0)=e^{iw_k(\alpha_0)}=e^{3iw_k(\alpha_0)}=(-1)^{{k\over2}+1}$. 
This proves $d)$. 

\smallskip
By Lemma 23, $b)$, the zeros of $\varphi'_k$ in $\D^+$ 
are ${1\over 2}+ic_1,\ldots,{1\over 2}+ic_n$, where 
$c_1,\ldots, c_n$ 
are some real numbers satisfying the inequalities
\begin{equation}\label{eq89}
c_1>b_1>c_2>b_2>\ldots>c_n>b_n.
\end{equation}
The function $\varphi'_k$ hence does not vanish at any point of 
$\rR_k^+$. Moreover the poles ${1\over 2}+ib_1,\ldots,{1\over 2}+ib_m$ of 
$\varphi_k$ are simple. It follows, by the local inversion theorem, 
that if $\tau$ is a point of $\rR_k^+$, the germ at $\tau$ of the set 
$\varphi_k^{-1}(\bP_1(\R))$  is the germ at $\tau$
of some real analytic curve. 
Furthermore, this curve is orthogonal to the vertical line with abscissa 
${1\over 2}$ at the point $\tau$ if $\tau$ is one of the points 
${1\over 2}+ib_j$, with $1\le j \le n$, and is orthogonal to the 
unit circle at the point $\tau$, if $\tau$ is one of the points 
$u_j$, with $0 \le j \le n$ : this follows from the fact that a biholomorphic map
is conformal, and from the fact that $\varphi_k$ has values in 
$({1\over2} + i \R)\cup\{\infty\}$ on the half-line ${1\over2} + i \R_+^*$ and values in 
the unit circle on the set of points of $\mathfrak{H}$ with modulus $1$ (Proposition 3).

Let us prove that the germ of $\rR_k^+$ at infinity in $\D^+$ 
is the germ of a real 
analytic curve, asymptotic to the vertical line with abscissa ${1\over 4}$. 
It is sufficient to prove that, for $t>0$ sufficiently large, 
there exists a unique real number $\sigma(t)\in[0,{1\over2}]$ 
such that $\sigma(t)+it \in \rR_k^+$, and that $\sigma(t)$
tends to ${1\over 4}$ when $t$ tends to $+\infty$. Note that we 
have the equivalences
\begin{equation}\label{eq90}
\varphi_k(z)\sim {i\B_k\over 4\pi}e^{-2\pi i z},
\end{equation}
\begin{equation}\label{eq91}
\varphi'_k(z)\sim {\B_k\over 2}e^{-2\pi i z},
\end{equation}
when $z \in\D^+$ approaches infinity. The relation $\eqref{eq90}$ 
implies that, for any $\varepsilon~\in~]0,{1\over8}]$, there exists 
$\T(\varepsilon)>0$ such that, 
for  $t \ge \T(\varepsilon)$, the function 
$$
x \mapsto \Im\varphi_k(x+~it)
$$ 
does not vanish in the intervals $[0,{1\over 4}-\varepsilon]$ 
and $[{1\over 4}+\varepsilon, {1\over2}]$ and has opposite signs in them,  
hence vanishes at at least one point of the interval 
$[{1\over 4}-\varepsilon,{1\over 4}+\varepsilon]$.
The relation $\eqref{eq91}$ implies that, for $t$ sufficiently large, 
the function 
$$
x\mapsto {d\ \over dx}\Im\varphi_k(x+it)
 = \Im{d\varphi_k(x+it)\over dx}=\Im\varphi'_k(x+it)
 $$
takes non-zero values with constant sign on the interval 
$[{1\over 4},{3\over 4}]$, 
hence the function $x\mapsto \Im\varphi_k(x+it)$ 
vanishes at at most one point of this
interval. This completes the proof of our  assertion about the germ 
of $\rR_k^+$ at infinity.

Let $\Gamma$ be a connected component of $\rR_k^+$. We deduce from the 
previous two paragraphs that $\Gamma$ is a real analytic curve, 
that $\varphi_k$ induces an analytic isomorphism from $\Gamma$ 
onto some connected subset of $\bP_1(\R)$. This subset is disjoint from
$]-1, 1[$ by Proposition 4 and its boundary points can only
be $1, -1$ or $\infty$. Hence either it is the canonical image
in $\bP_1(\R)$ of $[1, + \infty]$ or $[-\infty, -1]$, in which
case $\Gamma$ joins one of the 
points $u_0,\ldots,u_n$ to one of the points ${1\over2}+ ib_1,\ldots,{1\over 2}+ib_n$; 
or it is equal to $[1,+\infty[$ or $]-\infty,-1]$, in which case one of the 
points $u_0,\ldots,u_n$ is an endpoint of $\Gamma$, and $\Gamma$ 
has an infinite branch asymptotic to the vertical line with abscissa ${1\over 4}$.

It follows from the previous paragraph that the connected 
components $\Gamma_j$ of the points $u_j$, where 
\hbox{$0\le j\le n$}, are pairwise disjoint and that their union is
$\rR_k^+$. Since they do not intersect, necessarily $\Gamma_j$ joins 
$u_j$ to ${1\over2}+ i b_j$ for $1\le j\le n$ and $\Gamma_0$
has an infinite branch asymptotic to the vertical line with 
abscissa  ${1\over 4}$. We have thus proved assertions $e)$, $f)$, $g)$ and $h)$.

\bigskip
\noindent
{\bf Remark 18.}$-$ 
It follows from this proof that, for every $\lambda \in \R$, the roots of the
equation $\varphi(z) = \lambda$ in $\D$ are simple, and assertions $a)$
and $b)$ of Proposition~4 hold for $\lambda$.

\bigskip
\noindent
{\bf Remark 19.}$-$
A proof similar to that of Proposition 4 shows that, on each of the
$2n+3$ connected components of the interior of $\D - \rR_k$, 
$\varphi_k$ is injective.

\end{document}